\documentclass[12pt]{article}
\usepackage{marginnote}
\title{\vspace{-2 cm} Cairo  lattice with time-reversal non-invariant vertex couplings}
\usepackage{float}
\usepackage{array}
\usepackage{tabularx}
\usepackage{appendix}
\usepackage{makecell}
\usepackage{graphicx}
\usepackage{dcolumn}
\usepackage{bm}
\usepackage{amssymb}
\usepackage[T1]{fontenc}
\usepackage[utf8]{inputenc}
\usepackage[raggedright]{titlesec}
\usepackage{blindtext}
\usepackage{commath}
\usepackage{cite}
\usepackage{subcaption}
\usepackage{caption}
\usepackage[textwidth = 7in]{geometry}
\usepackage{mathtools}
\usepackage[dvipsnames]{xcolor}
\usepackage{enumitem}
\usepackage{amsmath}
\usepackage{amsthm}
\usepackage{geometry}
\usepackage{lipsum}
\usepackage{thmtools}
\usepackage{hyperref, cleveref}
\usepackage{setspace}
\titleformat{\paragraph}[hang]{\normalfont\normalsize\bfseries}{\theparagraph}{1em}{}
\titlespacing*{\paragraph}{0pt}{3.25ex plus 1ex minus .2ex}{0.5em}
\titleformat{\subparagraph}[hang]{\normalfont\normalsize\bfseries}{\thesubparagraph}{1em}{}
\titlespacing*{\subparagraph}{0pt}{3.25ex plus 1ex minus .2ex}{0.5em}
\makeatletter
\DeclareRobustCommand{\change}{%
	\@bsphack
	\leavevmode
	\color{magenta}%
	\@esphack
}
\DeclareRobustCommand{\stopchange}{%
	\@bsphack
	\normalcolor
	\@esphack
}
\makeatother

\newcommand{\e}{\mathrm{e}}

\stepcounter{secnumdepth}
\stepcounter{tocdepth}

\renewcommand{\thesection}{\arabic{section}}
\renewcommand{\theequation}{\arabic{equation}}

\DeclareMathOperator\arccosh{arccosh}

\pagebreak

\newtheorem{theorem}{Theorem}[section]

\renewcommand{\thesection}{\arabic{section}}
\renewcommand{\theequation}{\arabic{equation}}
\pagebreak

\author{Marzieh Baradaran$^{1}$\; and Pavel Exner$^{2,3}$}

\date{\small 1) Department of Physics, Faculty of Science, University of Hradec Kr\'alov\'e, Rokitansk\'eho 62, 500 03 Hradec Kr\'alov\'e, Czechia \\
	2) Doppler Institute for Mathematical Physics and Applied Mathematics, Czech Technical University, B\v rehov\'a 7, 11519 Prague, Czechia \\
	3) Department of Theoretical Physics, Nuclear Physics Institute, Czech Academy of Sciences, 25068 \v{R}e\v{z} near Prague, Czechia  \\
	\emph{marzie.baradaran@yahoo.com, exner@ujf.cas.cz}}

\begin{document}

\captionsetup[figure]{labelfont={bf},labelformat={default},labelsep=period,name={Fig.},font={normalsize}}
\captionsetup[table]{labelfont={bf},labelformat={default},labelsep=period,name={Table}}
\maketitle

\begin{abstract}
We analyze the spectrum of a periodic quantum graph of the Cairo lattice form. The used vertex coupling violates the time reversal invariance and its high-energy behavior depends on the vertex degree parity; in the considered example both odd and even parities are involved. The presence of the former implies that the spectrum is dominated by gaps. In addition, we discuss two modifications of the model in which this is not the case, the zero limit of the length parameter in the coupling, and the sign switch of the coupling matrix at the vertices of degree three; while different they both yield the same probability that a randomly chosen positive energy lies in the spectrum.

\end{abstract}

\section{Introduction}
\label{sect:intro}

Spectral properties of Laplacians on metric graphs, usually referred to as \emph{quantum graphs}, offer a number of mathematically interesting questions as well as many important applications \cite{BK13, KN22}.  The graph edges being identified with finite or semi-infinite intervals, the operator acts as $\psi_{j}\mapsto -\psi_{j}''$ on the $j$th edge. To make it a self-adjoint operator on the Hilbert space which is the orthogonal sum of $L^2$ spaces on the graph edges, one has to impose appropriate boundary conditions at the graph vertices. Those are by far not unique: for a given vertex $v$ in which $d_v$ number of edges meet, the self-adjointness is ensured provided the functions at the vertex are matched through the condition
\begin{equation}\label{genbc}
	(U_v-I)\psi(v)+i\ell(U_v+I)\psi'(v)=0,
\end{equation}
in which $\psi(v)$ and $\psi'(v)$ are respectively the vectors of the boundary values of the functions and their `outward' derivatives at the vertex, $U_v$ is a $d_v\times d_v$ unitary matrix, $I$ is the identity matrix, and 
$\ell>0$ is the parameter fixing the length scale. In the quantum graph literature most attention is paid to simple cases of such coupling, in particular, to those leading to continuity at the vertices. 

\smallskip

However, other types of vertex coupling may be also of interest, having in the first place in mind that they give rise to different physical models. This concerns, \emph{inter alia}, couplings that violate the time-reversal invariance. An example of such a coupling was proposed in \cite{ET18} motivated by a recent attempt to model the anomalous Hall effect in a quantum graph setting \cite{SK15, SV23}. This coupling has interesting properties, among them the dependence of the high-energy transport on vertex parity, with the consequences explored in \cite{ET18} and subsequent papers mentioned below. The examples discussed typically involved graphs, both finite and infinite, where all the vertices had the same degree. This motivates us to look at a situation where the degrees are different, combining different parities. Specifically, we are going to investigate an infinite periodic graph in the form of a \emph{Cairo lattice}, see~Fig.~\ref{CairoUnitCell} below, with the vertices of degrees three and four and two incommensurate edge lengths.

\smallskip

The band spectrum of such a lattice is dominated by gaps because any infinite path of the graph must pass through vertices of degree three which are opaque at high energies. This changes, however, if we send the length-scale parameter in such a coupling to zero changing it to the simplest matching condition usually called Kirchhoff. Another way to change the transport properties of the lattice is a small modification of the coupling of \cite{ET18} consisting of the sign change of the matrix $U_v$ in vertices of degree three. We find that while the couplings are different, they lead to similar overall transport properties manifested by the fact that the probability that a randomly chosen positive number lies in the spectrum, introduced in \cite{BB13} to characterize the so-called spectral universality, is the same in both cases. Note that the universality concept was introduced for graphs with Kirchhoff coupling but our present discussion shows it can apply in other situations too.

\smallskip

The paper is organized as follows. As a preliminary, we recall in the next section the necessary information about the vertex couplings. Then, in Sec.~\ref{sect:PosGenRR}, we introduce and discuss the Cairo lattice model with the coupling of \cite{ET18}, called here $R$; we find its spectrum and show how it changes in the Kirchhoff limit at the vertices of degree three. Sec.~\ref{sect:RminusR} is devoted to the modification mentioned, the replacement of the $R$ coupling by $-R$ at vertices of degree three. The results of each case are summarized at the end of each section. The functions appearing in the spectral conditions are quite complicated; to make the presentation lucid, we defer them to two appendices.

\section{The vertex coupling}
\label{s:vertex}

Let us first recall the vertex coupling introduced in \cite{ET18} and introduce its modification we will use. At a vertex $v$ of degree $d_v$ they are described by the condition \eqref{genbc} with the $d_v\times d_v$ unitary matrices $U_v$ of the circulant type,
\begin{equation}\label{umatrix}
	U_v =\pm R:= \begin{pmatrix}
		0 & \pm1 & 0 & \dots & 0 & 0\\
		0 & 0 & \pm1 & \dots & 0 & 0\\
		\vdots & \vdots &\vdots & \ddots & \vdots & \vdots\\
		0 & 0 & 0 & \dots & 0 & \pm1\\
		\pm1 & 0 & 0 & \dots & 0 & 0\\
	\end{pmatrix}\,,
\end{equation}
which are obviously non-invariant with respect to the transposition; hence the corresponding quantum graph dynamics violates the time-reversal invariance \cite{ET21}. Written in components, the condition \eqref{genbc} acquies the form
 \begin{equation}\label{couplings}
 	\pm(\psi_{j+1}\mp\psi_{j})+i\ell_{d_v}\,(\pm\psi_{j+1}^{\prime}+\psi_{j}^{\prime})=0,
   \end{equation}
where we have used for the sake of simplicity the symbols $\psi_j$ and  $\psi_{j}^{\prime}$,$\; j\in\mathbb{Z}\;(\textstyle{\rm mod}\;d_v)$ for the components of the vectors $\psi(v)$ and $\psi^{\prime}(v)$ of boundary values, respectively. Moreover, since our model, described in Sec.~\ref{sect:PosGenRR} below, contains vertices of different degrees, we have added the subscript $d_v\,$ to the length-scale parameter $\ell$ in \eqref{genbc}; specifically, we will consider the lengths $\ell_3$ and $\ell_4$ at the vertices of the appropriate degrees. For the sake of brevity, we refer in the following to the vertex conditions \eqref{couplings} with the upper and lower sign as to the $R$ and $-R$ coupling, respectively.

\smallskip

Let us first inspect properties of a single vertex considering a star-shaped graph of $d_v$ halflines supporting the Laplacian $H$ determined, as a self-adjoint operator, by the $\pm R$ couplings. Concerning the essential component of the spectrum, it is easy to check that it is absolutely continuous and coincides with the interval $[0,\infty)\,$; the quantity of interest in this part of the spectrum is the on-shell scattering matrix equal to $S(k) = \frac{k\ell_{d_v}-1 +(k\ell_{d_v}+1)U_v}{k\ell_{d_v}+1 +(k\ell_{d_v}-1)U_v}$ with $U_v$ given by \eqref{umatrix}, where the momentum variable $k$ is the square root of energy \cite{BK13}. The entries of $S(k)$ are found as in \cite{ET18} being
\begin{equation}\label{sij,onshell}
	S_{ij}^{ \pm R} (k) =\pm\frac{1-\eta^2}{1-\eta^{d_v}}\bigg\{ -\eta \, \frac{1-\eta^{d_v-2}}{1-\eta^2}\,\delta_{ij}+(1-\delta_{ij})\,\eta^{(j-i-1)(\textstyle{\rm mod}\;d_v)}          \bigg\},
\end{equation}
where $\eta:=\pm\tfrac{1-k\ell_{d_v}}{1+k\ell_{d_v}}$; note that the signs $\pm$ in the of $S_{ij}^{ \pm R} (k)$ correspond to those on the right-hand side of \eqref{sij,onshell} and that the variable $\eta$ sign depends on the choice of the coupling.

\smallskip

The scattering matrices corresponding to the ${\pm R}$ couplings behave differently in the high-energy asymptotic regime. It was observed in \cite{ET18} that in the $+R$ case the parity plays a decisive role: for an odd $d_v$ we have $\lim_{k\to\infty} S^R(k)=I$, implying that the particle coming to the vertex is fully reflected and its wave function remains thus confined to a single edge only. For an even $d_v$, on the other hand, the limit is not a multiple of $I$ and the scattering at high energies is non-trivial. Consequences of this fact have been explored in a number papers \cite{BEL23, BE22, EL20} and references therein.

\smallskip

The reason of the difference is not the parity itself, however, rather the absence or presence, respectively, of the Dirichlet component of the coupling, in the terminology of \cite{BK13}, corresponding to the eigenvalue $-1$ of $U_v$. Since eigenvalues of $R$ are the complex roots of unity, $-1$ is an eigenvalue if and only if $d_v$ is even. The consequences of its absence were illustrated in \cite{ET21} on the example of a coupling in which the matrix $R$ was modified by a phase factor. In this paper we consider another modification of that type, this time referring just to a change of sign. In this case there is no effective decoupling of the edges at high energies regardless of the vertex parity, because $-1$ is always an eigenvalue of $-R$, and consequently, $\lim_{k\to\infty} S^{-R}(k)$ is never a multiple of unit matrix.

\smallskip

The negative spectrum of $H$ is necessarily discrete and it can be easily checked that it is always nonempty for any $d_v\ge 3$. It is sufficient to use the Ansätze $\psi_{j}(x)=C_{j}\,\e^{-\kappa x}$ with $\kappa>0$ for the wave function components; plugging them into the conditions \eqref{couplings}, we get a system of $d_v$ linear equations the non-trivial solutions of which are determined by the characteristic equations
 \begin{equation}\label{NegStarG}
 	(-1-i \ell_{d_v}\kappa )^{d_v}-	(\mp 1\pm i \ell_{d_v}\kappa)^{d_v}=0,
 \end{equation}
where the upper and lower signs correspond to the $R$ and $-R$ couplings, respectively. It follows that the negative eigenvalues of a star graph Hamiltonian with the $\pm R$ couplings at the vertex equal
\begin{equation}\label{Neg,evRStarG}
	E_{\pm R,m}=-\kappa^2=-\ell_{d_v}^{-2}\,\Big(\tan^{2}\frac{m\pi}{d_v}\Big)^{\pm1},
\end{equation}
with $m$ running through $ 1,\cdots,[\tfrac{d_v}{2}]$ and $1,\cdots,[\tfrac{d_v-1}{2}]$ for odd and even $d_v\,$, respectively. Note that for even values of $d_v$, the sets $\{E_{\pm R,m}\}$ in \eqref{Neg,evRStarG} coincide, differing only in the numbering, which is also obvious from \eqref{NegStarG} given the fact that the two expressions in the second brackets differ only by sign.

\section{Cairo model with $R$ coupling at all vertices}
\label{sect:PosGenRR}

\begin{figure}[h]
 	\centering
 	\includegraphics[scale=0.50]{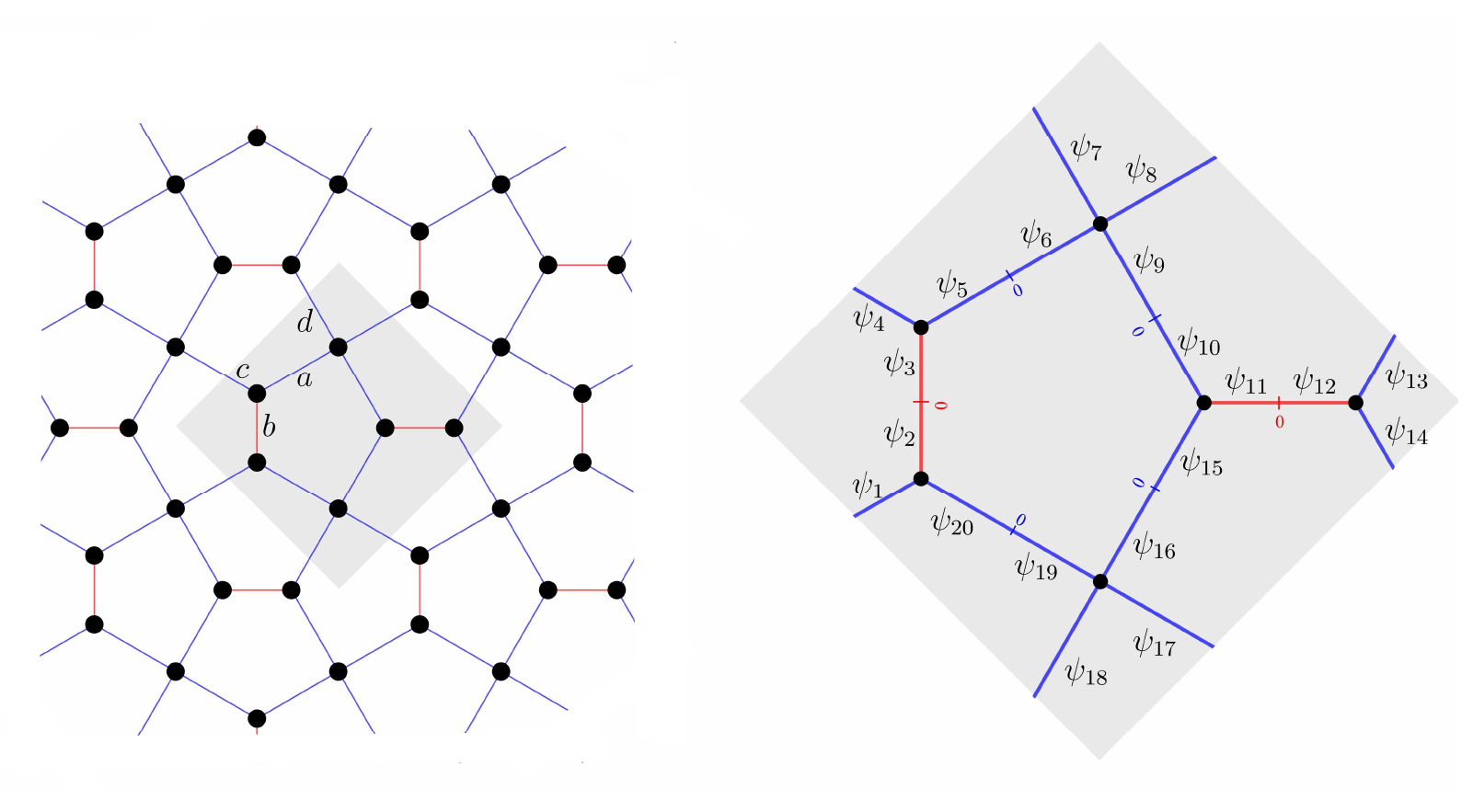}
 	\caption{Cairo lattice structure and its elementary cell indicated by the shaded area.}
 	\label{CairoUnitCell}
 \end{figure}

After the preliminary matters, let us proceed to describing the system of our interest. We consider an infinite quantum graph in the form of \emph{Cairo lattice} as sketched in Fig.~\ref{CairoUnitCell}. It is obtained by tiling the plane by non-equilateral pentagons with four edges of a length $a>0$ and one short one, $b=(\sqrt{3}-1)a$, marked respectively by blue and red colors in Fig.~\ref{CairoUnitCell}; accordingly, the long and short edges are (rationally) incommensurate. The system is periodic; an elementary cell of the model containing six vertices -- two vertices of degree four and four vertices of degree three -- is shown as the grey shaded area in Fig.~\ref{CairoUnitCell}. By elementary calculus, taking into account that the angles of the pentagons form the sequence $\frac{2\pi}3,\,\frac{2\pi}3,\,\frac{\pi}2,\,\frac{2\pi}3,\,\frac{\pi}2 $, one is able to find the lengths of the `loose-end' edges within the cell as $c=\frac a2(\sqrt{3}-1)$ and $d=\frac a2(3-\sqrt{3})$. Consequently, all the edge lengths in the elementary cell scale with the parameter $a$.
We begin with the `homogenous' situation where in all the vertices -- both of degree three and four -- the wave functions are matched through the $R$ condition. For the sake of clarity, we divide the following discussion of the spectral properties into several parts; the results are summarized in Theorem~\ref{thmGen} at the end of this section.

\subsection{Positive spectrum}	
	\label{sect:PosGen}
In view of the lattice periodicity, the spectral analysis can be reduced by the Floquet-Bloch decomposition theorem \cite[Chap.~4]{BK13} to investigation of an elementary cell of the graph. Since the motion is supposed to be free on the edges so that the Hamiltonian acts as $\psi(x)\mapsto -\psi''(x)$ on each component of the wavefunction, seeking positive energies $E=k^2>0\,$, we employ the Ans\"{a}tze $\psi_{\gamma}(x)=C_{\gamma}^{+}\e^{ikx} +C_{\gamma}^{-}\e^{-ikx}$ with $\gamma=1,2,...,20$ on the graph edges with the numbering shown in Fig.~\ref{CairoUnitCell}. To treat the negative spectrum, one has simply to replace $k$ by $i\kappa$ with $\kappa>0$.

\smallskip

As already mentioned, each edge of the graph is identified with an interval. We choose the coordinates on the graph edges to increase from `left to right' and from `bottom to top', thus the variable $x$ in $\psi _{1}(x)$ ranges through the interval $[-c,0]\,$, the variable in $\psi _{2}(x)$ and $\psi _{3}(x)$ ranges through $[-\frac b2,0]$ and $[0,\frac b2]$, respectively, etc. Obviously, this parametrization requires the functions to match smoothly at the midpoints of the edges in the interior of the elementary cell, in other words
\begin{equation}\label{midpoints}
	\psi _\gamma (0)=\psi _{\gamma+1}(0) \quad\text{and}\quad \psi _\gamma'(0)=\psi _{\gamma+1}'(0),
\end{equation}
for $\gamma=2,5,9,11,15,19$. In addition, since the model is periodic in two independent directions, the Floquet-Bloch decomposition imposes eight additional conditions at the loose ends of the edges, namely
\begin{equation}\label{Floq}
 	 \begin{aligned}
 	\psi _{8} ( d )&=\e^{i \theta_{1} } \psi _{1} (-c  ),\; & 	\psi _{8}' ( d )&=\e^{i \theta_{1} } \psi _{1}' (-c  ),
 	 \\
 	\psi _{13} ( c )&=\e^{i \theta_{1} } \psi _{18} (-d  ),\; & 	 \psi _{13}' ( c )&=\e^{i \theta_{1} } \psi _{18}' (-d  ),
 	 \\
 	\psi _{14} ( c )&=\e^{i \theta_{2} } \psi _{7} (-d  ),\; & 	\psi _{14}' ( c )&=\e^{i \theta_{2} } \psi _{7}' (-d  ),	
 	 \\
 	\psi _{17} ( d )&=\e^{i \theta_{2} } \psi _{4} (-c  ),\; & 	\psi _{17}' ( d )&=\e^{i \theta_{2} } \psi _{4}' (-c  )
 \end{aligned}
 \end{equation}
with the quasimomentum components running through the Brillouin zone, $\theta_1,\theta_2\in[-\pi,\pi)$. Furthermore, and most importantly, one has to impose the $R$ coupling condition \eqref{couplings} at all the six vertices of the elementary cell. For the sake of brevity, having in mind that the procedure is the same, we mention here the conditions only for two vertices of different parities. We choose the `clockwise' orientation at all the vertices. For instance, let us consider the vertex (of degree three) connecting the edges with the wavefunctions $\psi _{1}$, $\psi _{2}$ and $\psi _{20}$; imposing the $R$ coupling there and paying proper attention to the sign of derivatives that are taken in the `outward' direction, we get
 \begin{align*}\label{sys1}
&  	\psi _{2}(-\tfrac b2)-\psi _{1}(0)+i \ell_3 \left( \psi _{2}'(-\tfrac b2)-\psi _{1}'(0)\right) =0,\\
&  	\psi _{20}(-\tfrac a2)-\psi _{2}(-\tfrac b2)+i \ell_3 \left(  	\psi _{20}'(-\tfrac a2)+\psi _{2}'(-\tfrac b2)\right) =0,\\
&   \psi _{1}(0)-\psi _{20}(-\tfrac a2)+i \ell_3 \left( -\psi _{1}'(0)+\psi _{20}'(-\tfrac a2)\right) =0.
 \end{align*}
Next, let us consider the vertex (of degree four) at which the edges with the wave functions $\psi _{6}$, $\psi _{7}$, $\psi _{8}$ and $\psi _{9}$ meet; the matching conditions are
\begin{align*}
	&  	\psi _{7}(0)-\psi _{6}(\tfrac a2)-i \ell_4 \left( \psi _{7}'(0)+\psi _{6}'(\tfrac a2)\right) =0,\\
	&  	\psi _{8}(0)-\psi _{7}(0)+i \ell_4 \left(  	\psi _{8}'(0)-\psi _{7}'(0)\right) =0,\\
	&   \psi _{9}(-\tfrac a2)-\psi _{8}(0)+i \ell_4 \left( \psi _{9}'(-\tfrac a2)+\psi _{8}'(0)\right) =0,\\
	&   \psi _{6}(\tfrac a2)-\psi _{9}(-\tfrac a2)+i \ell_4 \left( -\psi _{6}'(\tfrac a2)+\psi _{9}'(-\tfrac a2)\right) =0.
\end{align*}
In a similar way, by imposing the $R$ coupling at the other four vertices of the elementary cell, we get thirteen more matching conditions. As a result, taking into account the constraints imposed by conditions \eqref{midpoints} and \eqref{Floq}, we arrive at a system of twenty linear equations which is solvable provided its determinant vanishes; this yields the spectral condition
\begin{equation}\label{spectGenPos}
f(k)=g(k) \left(\cos \theta_1 +\cos \theta_2 \right)+h(k) \left(\cos 2 \theta_1 +\cos 2 \theta_2 \right)+w(k) \,\cos \theta_1   \,\cos \theta_2 \,,
\end{equation}
where $f(k),\, g(k),\, h(k)$ and $w(k)$ are functions of the momentum variable $k$ depending also on the other parameters, $a$ and $\ell_{d_v}>0\,$, $d_v=3,4\,$.  Their explicit forms is cumbersome; they are given in Appendix~\ref{sect:appA} through Eqs. \eqref{f,app}-\eqref{w,app}. Inspecting the spectral condition \eqref{spectGenPos}, we infer that the spectrum consists of two types of spectral bands:

\begin{itemize}
\setlength{\itemsep}{-3pt}
\item The first are infinitely degenerate eigenvalues, \emph{flat bands} in the physicist's language, that appear when the spectral condition has solutions independent of the quasimomentum components $\theta_1$ and $\theta_2\,$. In the present case, it can happen only if $\ell_3=\ell_4$, then the number $k=\ell_3^{-1}$ belongs to the spectrum for $a=\frac{2 n\pi  }{\sqrt{3}+3}\ell_3$ with $n\in\mathbb{N}$. This can be checked by inspecting the functions $g(k),\, h(k)$ and $w(k)$ in \eqref{spectGenPos}, given by Eqs. \eqref{g,app}-\eqref{w,app}, checking whether they can vanish simultaneously. Consider first the simplest one, $h(k)$, which  obviously vanishes either for $k=\frac{ n\pi  }{a}$ or $k=\ell_3^{-1}$ provided $\ell_3=\ell_4$. Evaluating then $g(k)$, $w(k)$ and $f(k)$ at these values, we find that while the first two functions vanish, the third one, $f(k)$, reduces respectively to $8192  a^{-4}\pi ^4 \ell_3^4\cos ^4 \frac{\sqrt{3} \pi }{2} \neq0$ and $8192 \sin ^4 \frac{(\sqrt{3}+3) a}{2 \ell_3} $ which yields the claim.
\item Apart from the possible flat band mentioned above, the spectrum is absolutely continuous having a band-gap structure. To find it, we denote the right-hand side of \eqref{spectGenPos} by $R(\theta_1,\theta_2)$ and employ the Hessian matrix method for multivariable functions. First, we find that there exist six sets of critical points $p_i:=(\theta_1,\theta_2)\,$, $i=1,2,...,6$ -- note that $R(\theta_1,\theta_2)$ is symmetric with respect to the exchange of $\theta_1 $ and $\theta_2$. Denoting then, for the sake of brevity, $g:=g(k)$, $h:=h(k)$ and $w:=w(k)$, the critical points are as follows,
\begin{equation}\label{pointsExt}
	\begin{aligned}
	  p_1: \;\,&(0,\pm\pi),\\
	  p_2: \;\,&(0,0),\\
	  p_3: \;\,&(\pm\pi,\pm\pi),\\
	  p_4: \;\,&\left(\pm\pi,\arctan\Big(\frac{w-g}{h},\pm\frac{1}{h}\sqrt{(-g+4 h+w) (g+4 h-w)}\;\Big)\right),\\
	  p_5: \;\,&\left(0,\arctan\Big(-\frac{g+w}{h},\pm\frac{1}{h}\sqrt{(-g+4 h-w) (g+4 h+w)}\;\Big)\right),\\
	  p_6: \;\,&\bigg(\arctan\Big(\frac{-g}{4 h+w},\pm\sqrt\frac{(4 h+w)^2-g^2}{(4 h+w)^2}\;\Big),\arctan\Big(\frac{-g}{4 h+w}, \pm\sqrt\frac{(4 h+w)^2-g^2}{(4 h+w)^2}\;\Big)\bigg),
\end{aligned}
\end{equation}
where we have denoted $\arctan(x,y):=\arctan \frac yx$ taking into account which quadrant the point $(x,y)$ belongs to; the notation $(\pm,\pm)$ in \eqref{pointsExt} includes all four combination possibilities. Note also that the dispersion functions are even with respect to the quasimomentum components since $\theta_1$ and $\theta_2$ enter the arguments of cosine which is an even function. Computing next the determinant of the Hessian matrix, $\frac{\partial^2 R(\theta_1,\theta_2)}{\partial \theta_1^2}\,\frac{\partial^2 R(\theta_1,\theta_2)}{\partial \theta_2^2}- \frac{\partial^2 R(\theta_1,\theta_2)}{\partial \theta_1 \,\partial \theta_2}$, and inspecting the boundary of the Brillouin zone, we arrive at six local extrema of the function $R(\theta_1,\theta_2)$, namely
\begin{equation}\label{SIXextremas}
\begin{aligned}
\Lambda_1(k)&=2h-w ,\\
\Lambda_2(k)&=2g+2h+w ,\\
\Lambda_3(k)&=-2 g+2 h+w,\\
\Lambda_4(k)&= -\frac{g^2+8 g h-2 g w+w^2}{8 h} \qquad\text{\emph{iff}}\qquad (-g+4 h+w) (g+4 h-w)>0,\\
\Lambda_5(k)&= -\frac{g^2-8 g h+2 g w+w^2}{8 h} \qquad\text{\emph{iff}}\qquad (-g+4 h-w) (g+4 h+w)>0,\\
\Lambda_6(k)&=-\frac{g^2}{4 h+w}-2 h \qquad\qquad\qquad\text{\emph{iff}}\qquad (4 h+w)^2-g^2>0,
\end{aligned}
\end{equation}
in which each function $\Lambda_i(k)$ corresponds to the set of the critical points $p_i\,$, $i=1,2,...,6\,$, in \eqref{pointsExt}; needless to say, the constraints over the last three extrema are imposed by the real nature of the corresponding critical points. All of the above leads us thus to the conclusion that a number $k^2$ belongs to a spectral band if its square root, $k$, satisfies the conditions
\begin{equation}\label{BandGenPos1}
	\Lambda_i(k) \leq f(k)\leq    \Lambda_j(k),
\end{equation}
for some $i\neq j$ where $i, j \in \{1,2,...,6\}\,$, otherwise it belongs to a spectral gap. Let us emphasize one more time that one should pay proper attention to the accompanying constraints on the three extrema $\Lambda_i(k)$ with $i=4,5,6$ when using the band condition \eqref{BandGenPos1}. The band-gap pattern in dependence on the free parameter $a$ and the length scales $\ell_3$ and $\ell_4$ is illustrated in Figs.~\ref{Cairo,ell3=ell4=1a}--\ref{Cairo,ell3=1,a=1a}.
\end{itemize}

\subsection{High-energy asymptotics}
\label{sect:PosGenHigh}	
Let us now take a look at the spectrum in the high-energy regime, $k\to\infty$. Recall first that our Cairo model contains vertices of degree three and four at which the high-energy limit of the scattering matrix \eqref{sij,onshell} is respectively given by $\lim_{k\to\infty} S^{R}(k)=I_3$ and
\begin{equation} \label{splus4}
	\lim_{k\to\infty} S^{R}(k)=
\frac12 \left(\begin{array}{cccc}
	\phantom{-}1 &     \phantom{-}1 &       -1 &                  \phantom{-}1 \\
	\phantom{-}1 &     \phantom{-}1 &        \phantom{-}1 &     -1 \\
	-1 &               \phantom{-}1 &        \phantom{-}1 &       \phantom{-}1 \\
	\phantom{-}1 &     -1 &                  \phantom{-}1 &       \phantom{-}1 \\
\end{array}\right) .
\end{equation}
The latter limit indicates that the probabilities of leaving the vertex are asymptotically the same in any of the four directions, the former one obviously implies that the vertex becomes effectively decoupled at high energies. Consequently, as in the other examples of quantum graphs with the $R$ coupling at vertices of odd parities mentioned above, one expects the spectrum of the present model to be dominated by gaps which is indeed the case. To show that, we rewrite the spectral condition \eqref{spectGenPos} keeping the highest power of $k$, that is
\begin{equation} \label{asymptcond}
512 \;\ell_4^4\; \ell_3^8\,\left(\sin ^6 ka\;\cos ^2 ka\; \sin ^2 ( 1-\sqrt{3})ka \right)\, k^{12} +\mathcal{O}(k^{10})=0,
\end{equation}
from which we see that as $k\to\infty$ the spectral bands can appear only in the vicinity of the three numbers, $k=\frac{n\pi}{a}$, $(n-\frac12)\frac{\pi}{a}$ and $\frac{(\sqrt{3}+1)n\pi}{2a}$, $n\in\mathbb{N}$, at which the leading term in \eqref{asymptcond} vanishes. In particular, the probability that a randomly chosen positive energy lies in the spectrum, introduced by Band and Berkolaiko \cite{BB13} as
\begin{equation}\label{probsigma}
	P_{\sigma}(H):=\lim_{K\to\infty} \frac{1}{K}\left|\sigma(H)\cap[0,K]\right|,
\end{equation}
equals to zero in this case. Let us determine the asymptotic width of these bands.

\subsubsection*{(i) Narrow bands around the points $ \frac{n\pi}{a}\,$, $n\in\mathbb{N}$}

There may be six narrow bands in the vicinity of the roots of $\sin ^6 ka$, in general overlapping; to estimate their asymptotic width, we set $k_{1,n}=\frac{n\pi}{a}+\delta$, and consider the limit $n\rightarrow \infty\,$. Substituting $k_{1,n}$ into \eqref{spectGenPos}, and expanding the resulting equation around $\delta\to0$ to the leading order, we arrive at a sextic equation in $\delta$, in fact cubic for $\delta^2$, whose solutions are asymptotically of the form
$$   \delta_{n,j}=\pm\frac{\mathcal{F}_j(\tau_1;\ell_3,\ell_4,a)}{n}+\mathcal{O}(n^{-2}), \qquad j=1,2,3,$$
where $\tau_1:= (\cos\theta_1 +1)\,\cos\theta_2 +\cos\theta_1\in[-1,3]\,$ and the explicit form of the functions $\mathcal{F}_j$ is given in Appendix~\ref{sect:appA}. Since we have three pairs of solutions, each having both signs, each corresponding pair of bands has asymptotically the same width. Taking into account the range of the quasimomentume-dependent quantity $\tau_1$, the widths of the bands are
\begin{equation}
	\triangle E_{1,n} =
	\begin{cases}
		\frac{1}{3\,a\,\ell_3\,\ell_4} \,\Big(2 \sqrt{6 \ell_3 \sqrt{12 \ell_4^2+\ell_3^2} \cos \alpha+9\ell_4^2+6\ell_3^2}-6 \sqrt{\ell_4^2+\ell_3^2}\Big) +\mathcal{O}(n^{-1}) , \nonumber\\[10pt]
		\frac{1}{3\,a\,\ell_3\,\ell_4}\,\Big(2 \sqrt{3 \ell_3 \sqrt{12 \ell_4^2+\ell_3^2} \big(\sqrt{3} \sin \alpha-\cos\alpha\big)+9\ell_4^2+6\ell_3^2}-6 \sqrt{\ell_4^2+\ell_3^2}\Big) +\mathcal{O}(n^{-1}) , \nonumber\\[10pt]
		\frac{2}{3\,a\,\ell_3\,\ell_4}\,\sqrt{-3 \ell_3 \sqrt{12 \ell_4^2+\ell_3^2} \big(\sqrt{3} \sin\beta+\cos \beta\big)+9\ell_4^2+6\ell_3^2} +\mathcal{O}(n^{-1}), \nonumber
	\end{cases}
\end{equation}
in which $\alpha:=\frac{1}{3} \arctan \frac{6 \sqrt{3} \,\ell_4 \lvert \ell_3^2-4 \ell_4^2  \rvert }{36 \,\ell_4^2 \,\ell_3-\ell_3^3}$ and $\beta:=\frac{1}{3} \arctan \frac{6 \sqrt{3} \,\ell_4 \lvert \ell_3^3-4 \ell_4^2 \,\ell_3 \rvert }{36 \,\ell_4^2 \,\ell_3^2-\ell_3^4}$; as can be seen, the bands are of asymptotically constant width at the energy scale.

\subsubsection*{(ii) Narrow bands around the points $(n-\frac12)\frac{\pi}{a}\,$, $n\in\mathbb{N}$}
\label{sect:GenHigh2}

In this case, there are at most two narrow bands in the vicinity of the roots of $\cos ^2 ka$. As in the previous case, we set $k_{2,n}=(n-\frac12)\frac{\pi}{a} +\delta$ with $n\rightarrow \infty\,$; substituting this Ansatz into \eqref{spectGenPos} and keeping the leading term, we arrive at a quadratic equation for $\delta$ the solutions of which are of the form
$$ \delta_n =\pm \frac{1}{ n}\,\frac{\sqrt \tau_2}{2\pi\ell_3}+\mathcal{O}(n^{-2}) ,$$
where $\tau_2=(\cos\theta_2+1) (\cos\theta_1+1)\in[0,4]$. Consequently, both bands are again of a asymptotically constant width at the energy scale, namely
\begin{align}
	&  \triangle E_{2,n} = \frac{2}{  a\, \ell_3}+\mathcal{O}(n^{-1}). \nonumber
\end{align}
The question whether the gap between them is open depends on the error term; we leave it open.

\subsubsection*{(iii) Narrow bands around the points $\frac{(\sqrt 3+1)n\pi}{2a}\,$, $n\in\mathbb{N}$}

This time the bands appear in the vicinity of the roots of $\sin ^2 ( 1-\sqrt{3})ka$. In analogy with the previous cases, we set $k_{3,n}=\frac{(\sqrt 3+1)n\pi}{2a}+\delta$ with $n\rightarrow \infty\,$; substituting this into \eqref{spectGenPos} and keeping the leading term, we arrive at a quadratic equation for $\delta $ the solutions of which are now different for the even and odd values of $n$. For an even $n$, we have
$$   \delta_n=\frac{a\,(\sqrt{3}-1) }{ n^2\pi ^2 \,\ell_3^2}\left(2+4 \cos \sqrt{3} n\pi   -\tau_3^{\pm}\right)\,\,  \csc \sqrt{3} n\pi   +\mathcal{O}(n^{-3}),\qquad n\in2\mathbb{N},$$
in which $\tau_3^{\pm}:=  \cos \theta_1 +\cos \theta_2\,\pm\, 2 \sqrt{(\cos \theta_1 +1) (\cos \theta_2 +1)} \, $; these quantities, $\tau_3^{+}$ and $\tau_3^{-}$, run respectively through the intervals $[-2,6]$ and $[-2,0]$. Since these ranges overlap, the solutions merge into a single band of the width
\begin{equation}
	\triangle E_{3,n} =
		\frac{16}{ n \pi \,\ell_3^2}\, \lvert\csc \sqrt{3} n \pi\rvert+ \mathcal{O}(n^{-2}),\quad n\in2\mathbb{N}. \nonumber
\end{equation}
Similarly, for odd values of $n$ we get
$$   \delta_n=\frac{a\,(\sqrt{3}-1) }{ n^2\pi ^2 \,\ell_3^2}\left(-2+4 \cos \sqrt{3} n\pi-\tau_4\right)\,\,  \csc \sqrt{3} n\pi   +\mathcal{O}(n^{-3}),\quad n\in2\mathbb{N}-1,$$
where $\tau_4:=\cos \theta_1 +\cos \theta_2\in[-2,2]$; the bands width is then obtained as
$$\triangle E_{3,n} =	\frac{8}{ n \pi \,\ell_3^2}\, \lvert\csc \sqrt{3}  n \pi\rvert+ \mathcal{O}(n^{-2}),\quad n\in2\mathbb{N}-1.$$
The band widths are no longer asymptotically constant, but they remain \emph{bounded} as $n\to\infty$. Indeed, the $\csc$ modulus does not change if the argument is shifted by $\pi m$ with $m\in\mathbb{N}$. Using the fact that $\frac{\sin x}{x}\ge \frac12$ holds for small $|x|$ we get
$$\frac{1}{n}\,\big|\csc(\pi(\sqrt{3}n-m)\big| = \Big[n\sin\big(\pi n\big|\sqrt{3} -\textstyle{\frac{m}{n}}\big|\big)\Big]^{-1}
\le \Big[n\, \frac12\,\pi n\big|\sqrt{3} -\textstyle{\frac{m}{n}}\big|\Big]^{-1},  $$
however, $\sqrt{3}$ is an algebraic number of degree $2$, hence there is a $c>0$ such that $\big|\sqrt{3} -\textstyle{\frac{m}{n}}\big|\ge \frac{c}{n^2}$ holds for all coprime $m,n$, cf.~\cite{HW79, Sch91}.

\subsection{Negative spectrum}
\label{sect:NegGen}	

As already mentioned, the secular equation for the negative spectrum can be obtained by replacing $k$ by $i\kappa$ with $\kappa>0$ in the spectral condition \eqref{spectGenPos}. Mimicking the argument of Sec.~\ref{sect:PosGen}, we infer that a negative number $-\kappa^2$ belongs to a spectral band if $\kappa$ satisfies the conditions
\begin{equation}\label{BandGenNeg1}
	\Lambda_i(\kappa) \leq f(\kappa)\leq    \Lambda_j(\kappa),
\end{equation}
for some $i\neq j \in \{1,2,...,6\}$. The $\Lambda_i(\kappa)$ and $f(\kappa)$ here are obtained from \eqref{SIXextremas} and \eqref{f,app}, respectively, by replacement of the real $k$ by imaginary one which changes sine and cosine to their hyperbolic versions; needless to say, this is also the case for the other three functions $g(\kappa),\, h(\kappa)$ and $w(\kappa)$ included in $\Lambda_i(\kappa)$. As a result of \eqref{BandGenNeg1}, the negative spectrum, too, has a band-gap structure, this time consisting of a finite number of bands and gaps, cf. the first bullet point bellow. The band-gap pattern is illustrated in Figs.~\ref{Cairo,ell3=ell4=1a}--\ref{Cairo,ell3=1,a=1a}. We see that the negative spectrum has the following properties:

\begin{itemize}
\setlength{\itemsep}{-3pt}
\item Concerning the number of negative bands, we have shown in Sec.~\ref{s:vertex} that the negative eigenvalues of a star graph with $d_v\ge 3$ semi-infinite edges connected through $R$ coupling are given by \eqref{Neg,evRStarG}, in particular, we have $E_{+ R}=-3\ell_3^{-2}$ and $-\ell_4^{-2}$ for  $d_v=3$ and $4$, respectively; recall the general result by which the number of negative eigenvalues of a star graph coincides with the number of eigenvalues of $U_v$ in the upper complex halfplane \cite[Thm.~2.6]{BET22}. The lattice spectral bands shrink to these eigenvalues in the limit $a\to\infty$ when the tunelling between the vertices becomes negligible as it seen in Figs.~\ref{Cairo,ell3=ell4=1a} and \ref{Cairo,ell3=2,ell4=23a}. Different bands can shrink to the same eigenvalue, of course, but given the number of vertices in the elementary cell and the dimension of the corresponding resolvent singularities, we can have at most six negative spectral bands; in reality there are at most three.
\item As it can be expected, the band shrinking to the eigenvalues $E_{+R}=-\ell_4^{-2}$ and $-3\ell_3^{-2}$ in the limit $a\to\infty$ is exponentially fast. To see that, we note that $\cosh (\pm x)\approx \tfrac12\,e^x$ and $\sinh  (\pm x)\approx \pm\tfrac12\,e^x$ holds as $x\to\infty$ and rewrite the spectral condition as
\begin{equation}\label{Neg,gen,lead}
\frac{1}{2} \left(1-\kappa ^2 \ell_4^2\right)^2 \left(3-\kappa ^2 \ell_3^2\right)^4 \, \e^  {2 (\sqrt{3}+3) \kappa a}+\mathcal{O}(\e^  {8 \kappa a})=0,
\end{equation}
which can be satisfied only if the first expression on the left-hand side matches the error term. To find the approximate band width, we set $\kappa_1 =\ell_4^{-1}+\delta $; expanding then the resulting equation in the vicinity of $\delta=0$, we get in the leading order a linear equation for $ \delta $ the solution of which yields
$$   E_1= -\frac1{\ell_4^{2}}-\frac2{\ell_4^2}\,\frac{ \,\tau_1\,\ell_4^2 - \ell_3^2}{3 \ell_4^2-\ell_3^2}\,\e^{-\frac{2 a}{\ell_4}}+\mathcal{O}\big(\e^{\frac{(\sqrt{3}-3) a}{\ell_4}}\big),$$
where $\tau_1= (\cos\theta_1 +1)\,\cos\theta_2 +\cos\theta_1\,\in[-1,3]\,$, and as a result, the band width is given by $\triangle E_1=\frac{8 }{3 \ell_4^2-\ell_3^2}\,\e^{-\frac{2 a}{\ell_4}}+\mathcal{O}\big(\e^{\frac{(\sqrt{3}-3) a}{\ell_4}}\big)$. Putting similarly $\kappa_2 =\sqrt{3}\,\ell_3^{-1}+\delta $, we get from the spectral condition \eqref{Neg,gen,lead} a quadratic equation in $\delta$ the solutions of which yield
\begin{align*}
	 E_2&= -\frac3{\ell_3^{2}}\mp\,\frac{2 \sqrt{2}}{\ell_3^{2}}\,\e^{\frac{(\sqrt{3}-3) a}{\ell_3}}-\frac{4}{3 \sqrt{3} \ell_3^3}\Big(3  (\sqrt{3}-5 ) a-\frac{3 \sqrt{3} \ell_3  (\ell_3^2-5 \ell_4^2 )}{\ell_3^2-3 \ell_4^2} \Big)\,\e^{\frac{2(\sqrt{3}-3) a}{\ell_3}}\\
	 &\pm\,\frac{2 \sqrt{2}  \left(3 \ell_4^2+\ell_3^2\right)}{\ell_3^4-3 \ell_4^2 \ell_3^2}\,(\cos\theta_1 +\cos\theta_2)\,\e^{\frac{-2\sqrt{3} a}{\ell_3}}+\mathcal{O}\big(\e^{\frac{- (\sqrt{3}+3) a}{\ell_3}}\big)
\end{align*}   	
with the first two terms independent of $\theta_1,\theta_2$, and as $-2\leq\cos\theta_1 +\cos\theta_2\leq2$, the width of the bands shrinking to $-3\ell_3^{-2}$ is  $\triangle E_2=\frac{8 \sqrt{2}  \left(3 \ell_4^2+\ell_3^2\right)}{\ell_3^4-3 \ell_4^2 \ell_3^2}\,\e^{\frac{-2\sqrt{3} a}{\ell_3}}+\mathcal{O}\big(\e^{\frac{- (\sqrt{3}+3) a}{\ell_3}}\big)$. Note that the error terms can be chosen in both cases independent of the quasimomentum.
\item There are \emph{no flat bands} in the negative part of the spectrum. To see that, one has to check that the quasimomentum-dependent terms in the spectral condition, the functions $g(\kappa),\, h(\kappa)$ and $w(\kappa)$ cannot vanish simultaneously. It is sufficient to inspect the simplest one, $h(\kappa)\,$. By a simple manipulation, the condition $h(\kappa)=0$ can be brought to the form
$$ -\left(\frac{\kappa  (\ell_4-\ell_3) }{ \kappa ^2 \,\ell_3\,\ell_4+1 }\right)^2=\tanh ^2  \kappa a, $$
which cannot hold even if $\ell_3=\ell_4$, because the right-hand side is positive for any $\kappa>0$.
\end{itemize}

\begin{figure}[H]
	\centering
	\subfloat[   the model with the $R$ coupling  ]{{\includegraphics[scale=0.55]{ 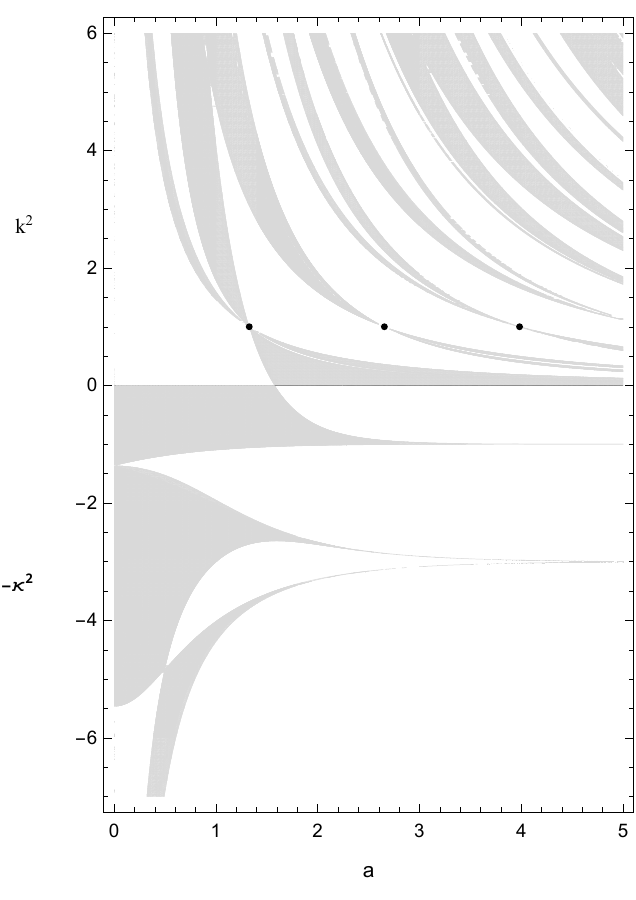  }} \label{Cairo,ell3=ell4=1a}}
	\hspace{50pt} 
	\subfloat[  the model with the $(-1)^{d_v}R$ coupling ]{\includegraphics[scale=0.55]{  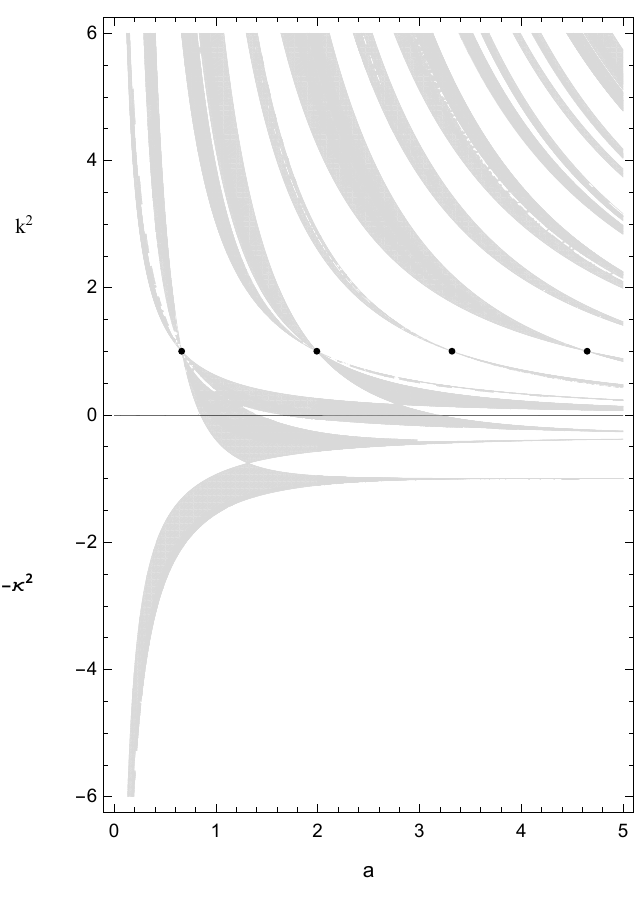  }\label{Cairo,ell3=ell4=1b}}
	\caption{  Spectrum of the Cairo  lattice in dependence on the parameter $a$ with  $\ell_3=\ell_4=1$. In this and the subsequent figures~\ref{Cairo,ell3=2,ell4=23}--\ref{Cairo,ell3=1,a=1}, the left picture corresponds to the model in which $R$ coupling is imposed at all vertices (discussed in Sec.~\ref{sect:PosGenRR}) while the picture on the right corresponds to the model where $R$ coupling is imposed at vertices of degree four, and $-R$ coupling is at vertices of degree three (discussed in Sec.~\ref{sect:RminusR}). The black dots denote the flat bands.}
	\label{Cairo,ell3=ell4=1}
\end{figure}


\begin{figure}[H]
	\centering
	\subfloat[    the model with the $R$ coupling  ]{{\includegraphics[scale=0.37]{ 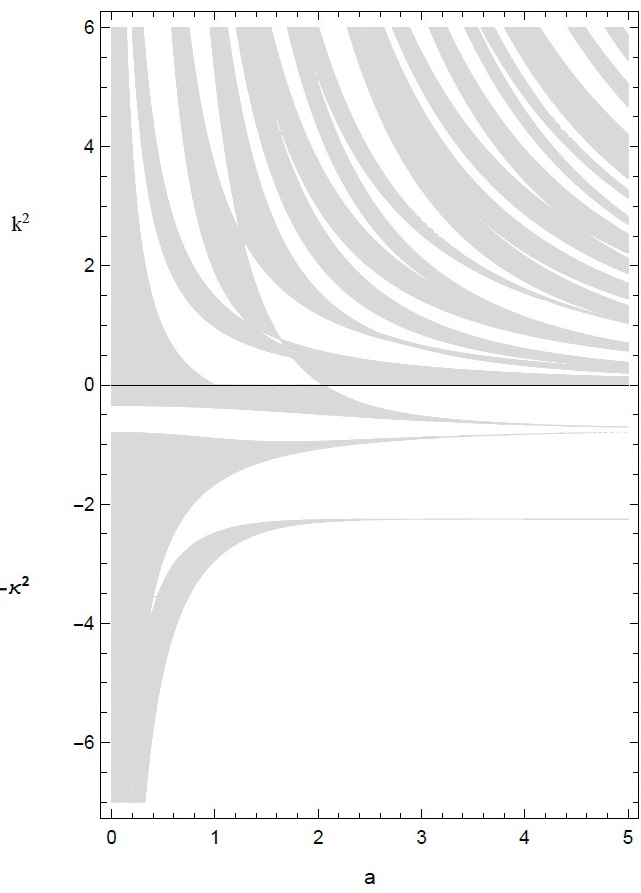  }} \label{Cairo,ell3=2,ell4=23a}}
	\hspace{50pt} 
	\subfloat[   the model with the $(-1)^{d_v}R$ coupling ]{\includegraphics[scale=0.37]{ 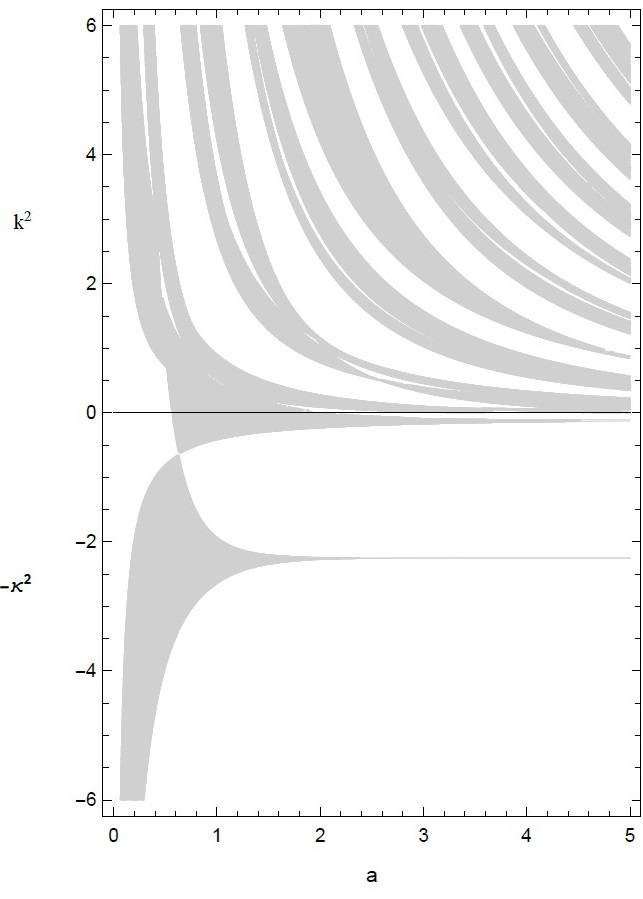   }\label{Cairo,ell3=2,ell4=23b}}
	\caption{   Spectrum of the Cairo  lattice in dependence on the parameter $a$ with  $\ell_3=2$ and $\ell_4=\frac23$.}
	\label{Cairo,ell3=2,ell4=23}
\end{figure}


\begin{figure}[H]
	\centering
	\subfloat[  the model with the $R$ coupling ]{{\includegraphics[scale=0.37]{ 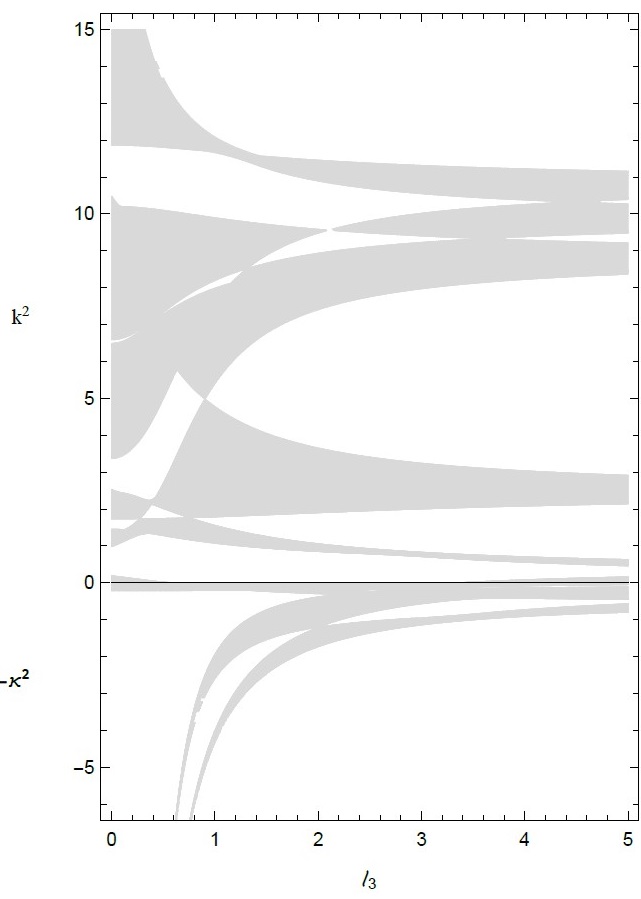  }} \label{Cairo,ell4=2,a=1a}}
	\hspace {50pt} 
	\subfloat[  the model with the $(-1)^{d_v}R$ coupling ]{\includegraphics[scale=0.37]{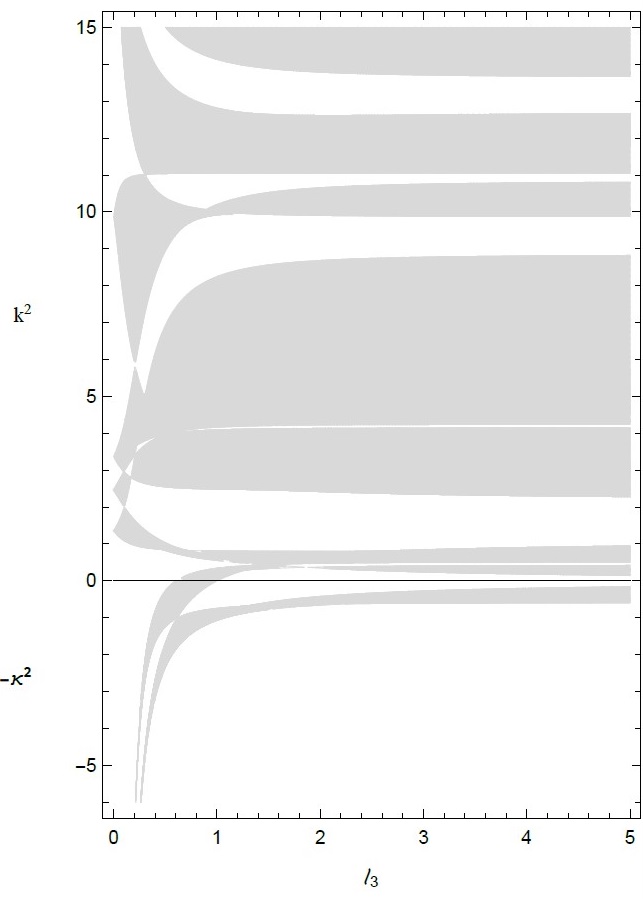  }\label{Cairo,ell4=2,a=1b}}
	\caption{   Spectrum of the Cairo  lattice in dependence on the length scale $\ell_3$ with $\ell_4=2$ and $a=1$.}
	\label{Cairo,ell4=2,a=1}
\end{figure}


\begin{figure}[H]
	\centering
	\subfloat[  the model with the $R$ coupling  ]{{\includegraphics[scale=0.4]{ 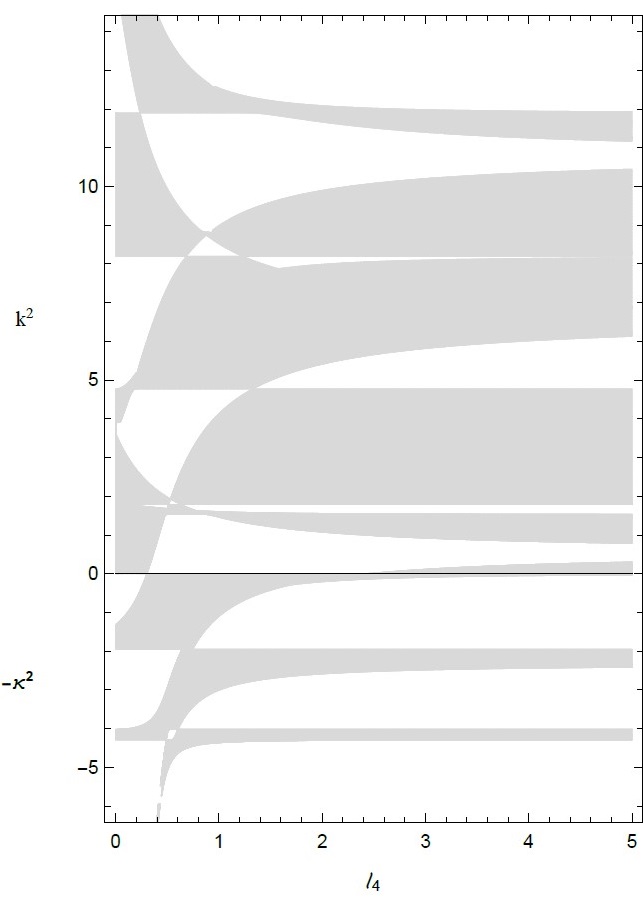  }} \label{Cairo,ell3=1,a=1a}}
	\hspace{50pt} 
	\subfloat[  the model with the $(-1)^{d_v}R$ coupling ]{\includegraphics[scale=0.4]{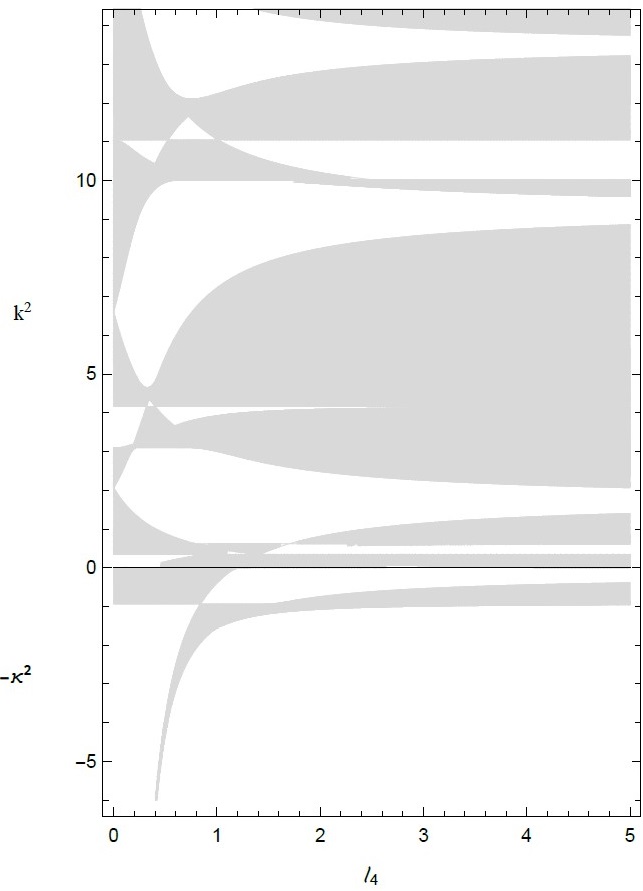  }\label{Cairo,ell3=1,a=1b}}
	\caption{   Spectrum of the Cairo  lattice in dependence on the length scale $\ell_4$ with $\ell_3=1$ and $a=1$. }
	\label{Cairo,ell3=1,a=1}
\end{figure}

\subsection{The Kirchhoff limit, $\ell_3\to 0$}
\label{sect:Ell30}

The spectral condition \eqref{spectGenPos} depends continuously on the parameters of the model, even analytically \cite[Thm.~2.5.4]{BK13}, and so do in view of implicit function theorem the eigenvalues of the fiber operators in the Floquet-Bloch decomposition. The limits of such functions are of particular interest. One often considers edges shrinking to a point \cite{BLS19, Bo21}; here in contrast we will look what happens if the length parameter $\ell_3$ in \eqref{couplings} with the upper signs tends to zero. Summing this condition over $j$, we get $\sum_j \psi'_j = 0$, and the limit $\ell_3\to 0$ yields the continuity, $\psi_j=\psi_{j+1}$, i.e. we arrive at what is usually called the Kirchhoff coupling. It is clear that the limit changes the spectral character substantially, a single vertex has no longer a negative spectrum which in accordance with \eqref{Neg,evRStarG} escapes to $-\infty$, while in the positive part the scattering is independent of $k$ and the edges are no longer asymptotically decoupled at high energies. Let us look what the limit does with the Cairo lattice spectrum.

\subsubsection{Positive spectrum}	
\label{sect:PosEll30High}

We take the limit $\ell_3\to 0$ in the spectral condition \eqref{spectGenPos}, that is, in the functions $f(k)$, $g(k)$, $h(k)$ and $w(k)$ in Eqs.~\eqref{f,app}-\eqref{w,app}; naturally, the condition of belonging to the spectral bands is still given by \eqref{BandGenPos1}. The band-gap pattern changes as illustrated in Fig.~\ref{Cairo,ell3=0,ell4=1}, in particular, the spectral bands are now wider. To get this claim a proper meaning, let us take a closer look at the spectrum behavior in the high-energy regime and the probability \eqref{probsigma} of belonging to the spectral bands; recall that for any $\ell_3>0$ this quantity was equal to zero.

\begin{figure}[H]
		\centering
	\subfloat[ The region in which the condition \eqref{con,ell30,xy1} with the functions \eqref{fghtilde,xy} holds; this corresponds to the Cairo model with the $R$ coupling in the limit $\ell_3\to0$.] {{\includegraphics[width=0.40\hsize]{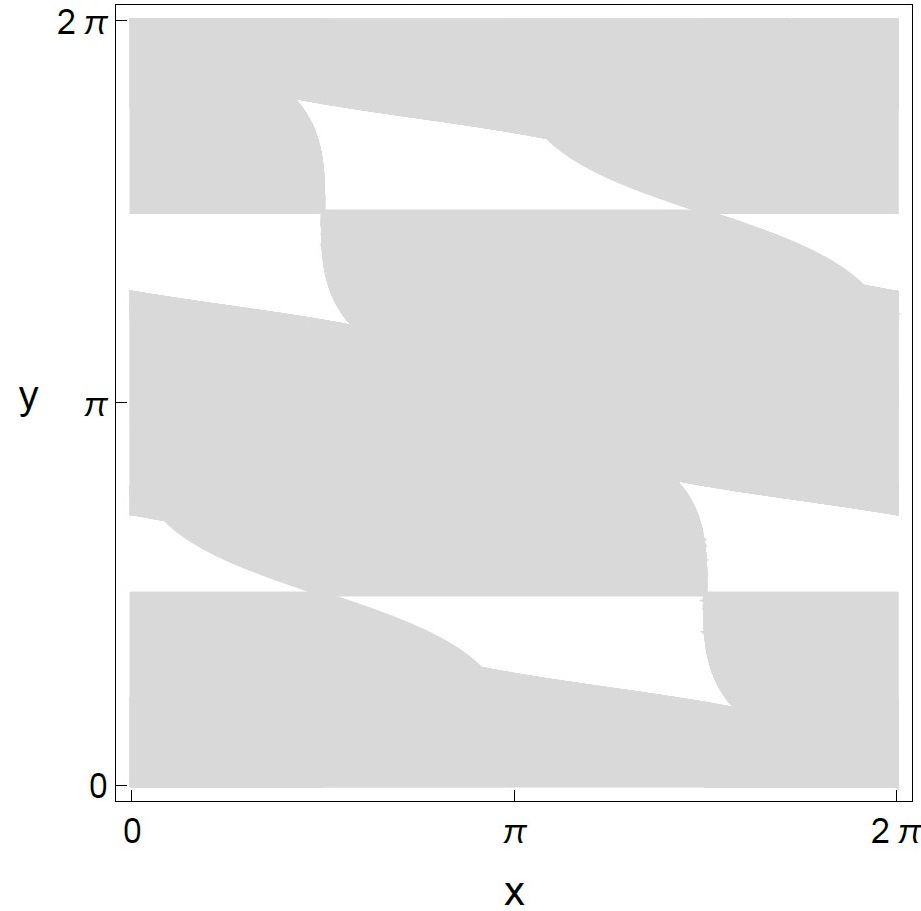}} \label{FigProbRR}}
	\hspace{70pt} 
	\subfloat[The region in which the condition \eqref{con,ell30,xy1} with the functions \eqref{fghtilde,xy,minusR} holds; this corresponds to the Cairo model with the $(-1)^{d_v}R$ coupling.]{\includegraphics[width=0.40\hsize]{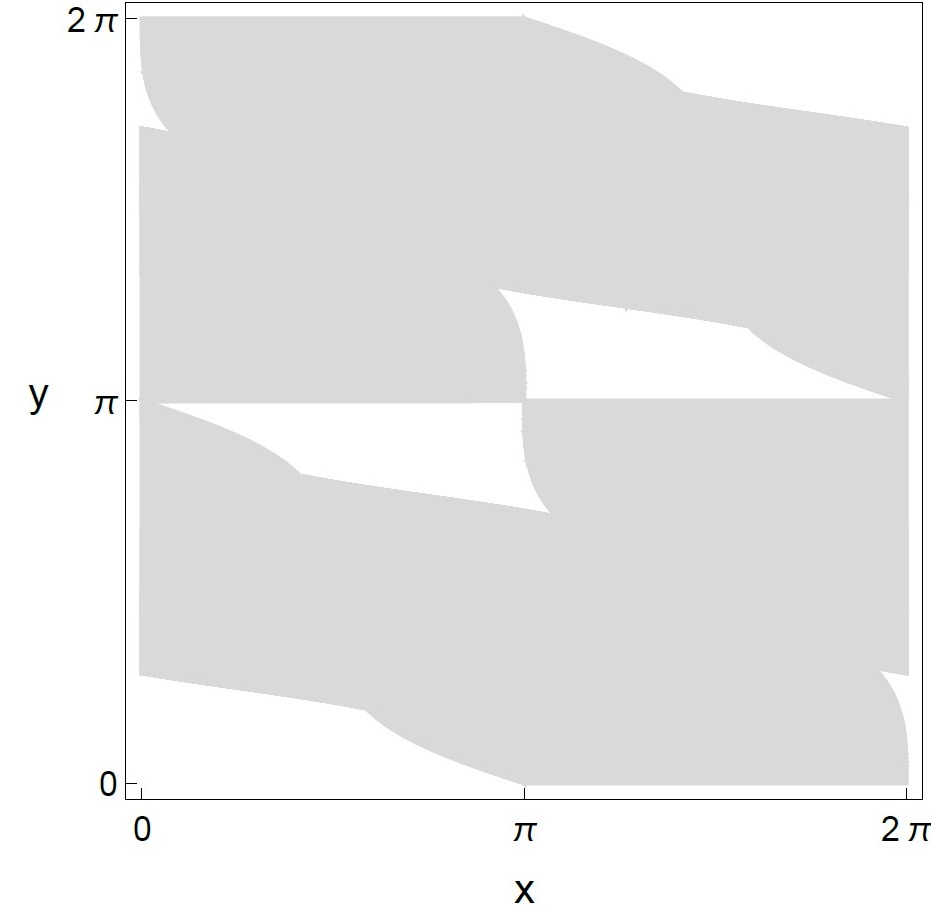}\label{FigProbRminusR}}
	\caption{The grey shaded area equals to $4\pi^2\,P_{\sigma}(H)$ determined by the conditions \eqref{con,ell30,xy1}; the axes correspond to $x:=\sqrt{3}\, ka$ and $y:=ka$ in the regime $k\to\infty$.}
	\label{FigProbsArea}
\end{figure}

Keeping the highest power of $k$ in the spectral condition as $k\to\infty$, we rewrite it in the form
\begin{align} \label{asymptcond,ell30}
	4 \;\ell_4^4\; \cos ^2 ka\,&\big\{	\tilde{f}(k)-\tilde{g}(k) (\cos \theta_1 +\cos \theta_2 )-\tilde{h}(k) (\cos 2 \theta_1 +\cos 2 \theta_2 )-\tilde{w}(k) \,\cos \theta_1   \,\cos \theta_2 \big\}\, k^{4} \\ \nonumber
	&+\mathcal{O}(k^{2})=0,
\end{align}
in which
\small
\begin{equation}\label{fghwtilde}
\begin{aligned}
\tilde{f}(k)&=-54 \cos 2ka+60 \cos 4ka-18 \cos6ka+\cos  2 (\sqrt{3}-4)ka +27 \cos  2 \sqrt{3}ka \\
& +16 \cos  (\sqrt{3}-3) ka -6 \cos 2 (\sqrt{3}-3)ka+11 \cos 2 (\sqrt{3}-2)ka-32 \cos (\sqrt{3}-1) ka \\
& -6 \left(2 \cos 2  (\sqrt{3}-1 ) ka+8 \cos  (\sqrt{3}+1 ) ka+9 \cos 2 (\sqrt{3}+1)ka-10\right)+81 \cos  2 (\sqrt{3}+2) ka , \\[5pt]
\tilde{g}(k)&=- 8 \cos 2ka+8 \cos  (\sqrt{3}-5 ) ka-32 \cos   (\sqrt{3}-3 ) ka+36 \cos  2 \sqrt{3}ka +56 \\
& +4 \cos  2 (\sqrt{3}-2)ka+32 \cos (\sqrt{3}-1) ka-24 \cos  2 (\sqrt{3}-1) ka-72 \cos  (\sqrt{3}+3)ka ,\\[5pt]
\tilde{h}(k)&=-16 \cos ^2 ka, \\[5pt]
\tilde{w}(k)&=-40 \cos 2ka-16 \cos  (\sqrt{3}-3)ka+4 \cos  2 (\sqrt{3}-2) ka+36 \cos  2 \sqrt{3} ka+24 \\
& +32 \cos (\sqrt{3}-1) ka-24 \cos 2 (\sqrt{3}-1)ka+48 \cos (\sqrt{3}+1) ka.
\end{aligned}
\end{equation}
\normalsize
The bands naturally appear where the leading term in \eqref{asymptcond,ell30} matches the error one, that is, around the points where the leading term vanishes. Consequently, there are two types of spectral bands. The first are narrow bands in the vicinity of the roots of $\cos ^2 ka\,$. We can specify their asymptotic behavior as in Sec.~\ref{sect:PosGenHigh} but we will not do that. What is more important, are the wide bands around the points determined by vanishing of the expression in the curly bracket in \eqref{asymptcond,ell30}. We are not that much interested in the bands themselves, rather in the global quantity $P_\sigma(H)$. Instead of trying to derive it from \eqref{BandGenPos1} directly, we use the spectral universality observed in \cite{BB13}. This is made possible by the incommensurability of the graph edges, $a$ and $b$, which allows us to regard the quantities $\sqrt{3}\, ka=:x$ and $ka=:y$ in the arguments of trigonometric functions in \eqref{fghwtilde} as a pair of independent identically distributed random variables in $[0,2\pi)$. In these new variables, the functions entering \eqref{fghwtilde} become
\small
\begin{equation}\label{fghtilde,xy}
	\begin{aligned}
		\tilde{f}(x,y) &=  32 \, \left(\cos (x-2 y)-3 \cos x\right)\cos y-6 \cos 2 (x-3 y)+24 \left(5-3 \cos 2y\right) \cos ^2 2 y  \\
		&   + \left(92 \cos 4y-66 \cos 2y+27\right)\cos 2x+14\left(3 \sin 2y-5 \sin 4y\right)\sin 2x +\cos 2 (x-4 y)  , \\[5pt]
		\tilde{g}(x,y) &=-24 \cos 2(x-y)-32 \cos (x-3 y)+4 \cos 2(x-2 y)+8 \cos (x-5 y)+36 \cos 2x+56\\
		& +32 \cos (x-y)-72 \cos (x+3 y)-8 \cos 2y , \\[5pt]
		\tilde{h}(x,y) &=-16 \cos ^2 y , \\[5pt]
		\tilde{w}(x,y) &=80 \cos x \,\cos y-24 \cos2 ( x- y)-16 \cos (x-3 y)+4 \cos 2( x-2 y)+36 \cos 2x+24\\
		&	-16 \sin x\,\sin y-40 \cos 2y.
	\end{aligned}
\end{equation}
\normalsize
Mimicking now the argument used in Sec.~\ref{sect:PosGen}, we infer that points of the spectrum have to satisfy the condition
\begin{equation}\label{con,ell30,xy1}
	\tilde{\Lambda}_i(x,y) \leq \tilde{f}(x,y)\leq    \tilde{\Lambda}_j(x,y),
\end{equation}
in which $x,y\in[0,2\pi)$ and $i\neq j \in \{1,2,...,6\}$; the functions $\tilde{\Lambda}_i(x,y)$ are nothing but $\Lambda_i(k)$ of \eqref{SIXextremas} in which $g$, $h$, and $w$ have been respectively replaced by $\tilde{g}(x,y)$, $\tilde{h}(x,y)$, and $\tilde{w}(x,y)$ from \eqref{fghtilde,xy}. Since $x$ and $y$ are identically distributed, the sought probability \eqref{probsigma} is nothing but the relative area of the region specified by conditions \eqref{con,ell30,xy1} in the square $(0,2\pi)^2$. The region is shown in Fig.~\ref{FigProbRR} being indicated by gray color; its area can be computed numerically, and being divided by $4\pi^2$ it yields the probability $P_{\sigma}(H)\approx 0.82$. This shows that while the spectral pattern converges locally as $\ell_3\to 0$, the limit is not uniform because $P_{\sigma}(H)=0$ holds for any $\ell_3>0$.

\subsubsection{Negative spectrum}	
\label{sect:NegEll30}

We obtain it by taking the limit $\ell_3\to 0$ in the spectral condition \eqref{BandGenNeg1}. There is a single negative spectral band as illustrated by Fig.~\ref{Cairo,ell3=0,ell4=1}, see also Fig.~\ref{Cairo,ell4=2,a=1a} where $\ell_4=2$. For large values of the parameter $a$, this band shrinks exponentially fast to the eigenvalue of a star graph of degree four given by~\eqref{Neg,evRStarG}, that is, $-\ell_4^{-2}$. To see that, we note that in the leading order of $a$, the spectral condition reads
 $$\frac{81}{2} \e^{2  (\sqrt{3}+3 ) \kappa a }  (\kappa ^2 \ell_4^2-1 )^2+\mathcal{O}(\e^{8 \kappa a })=0.$$
To find the asymptotic expression of the band, we set $ \kappa =\ell_4^{-1}+\delta $ and substitute it in the spectral condition. Repeating the argument of Sec.~\ref{sect:NegGen}, we find that the energy is
$$ E =-\frac1{\ell_4^2}-\frac2{3\ell_4^2}\,\tau_1\,\e^{-\frac{2 a}{\ell_4}} +\mathcal{O}\big(\e^{-\frac{(\sqrt{3}+1) a}{\ell_4}}\big) ,$$
where we have denoted again $\tau_1= (\cos\theta_1 +1)\,\cos\theta_2+\cos\theta_1 \in[-1,3]$; the band width asymptotics is thus $ \triangle E=\frac{8 }{3 \ell_4^2}\,\e^{-\frac{2 a}{\ell_4}} +\mathcal{O}\big(\e^{-\tfrac{(\sqrt{3}+1) a}{\ell_4}}\big).$


\begin{figure}[h]
	\centering
	\includegraphics[scale=0.40]{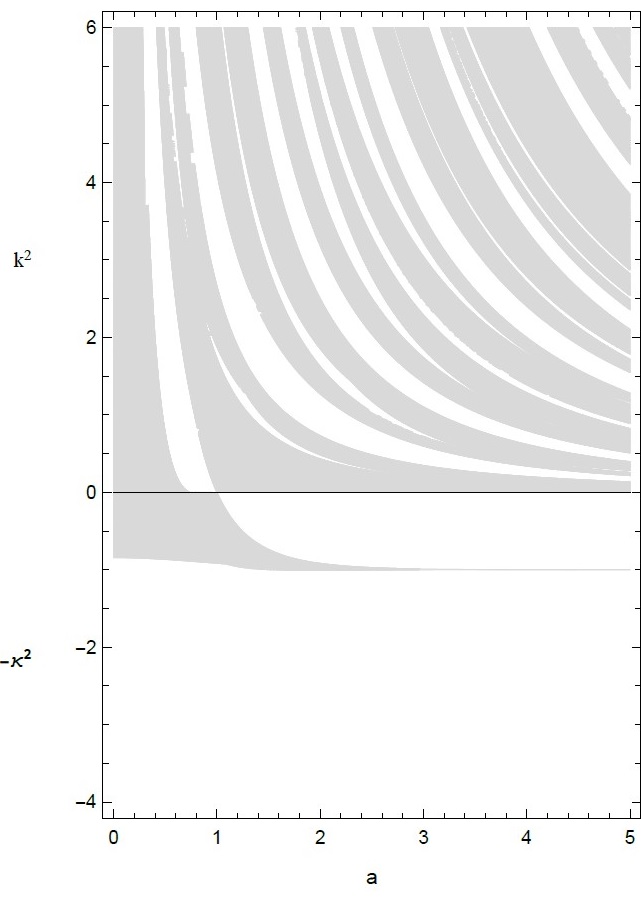}
	\caption{ Spectrum of the Cairo lattice with the $R$ coupling at all the vertices in dependence on the $a$ with $\ell_4=1$ in the limit $\ell_3\to0$.}
	\label{Cairo,ell3=0,ell4=1}
\end{figure}

\subsection{The results summary}
\label{sect:mainRR}

Before concluding this section, let us summarize our observations. Let $H^{+R}_{ \ell_3,\ell_4 }$ be the quantum graph Hamiltonian described in the beginning of this section; then, the obtained results can be summarized in the following theorem:
\begin{theorem}\label{thmGen}
The spectrum of $H^{+R}_{ \ell_3,\ell_4 }$ has the following properties:
	\begin{enumerate}[label=\textnormal{(\roman*)}]
		\setlength{\itemsep}{-3pt}
		\item \label{th1} The number $k=\ell_3^{-1}$, $\ell_3>0$ always belongs to $\sigma(H^{+R}_{\ell_3,\ell_3 })$ for $a=\frac{2 n\pi  }{\sqrt{3}+3}\ell_3$ with $n\in\mathbb{N}$.
		\item \label{th2} Away of the above infinitely degenerate eigenvalue, the spectrum is absolutely continuous having a band-and-gap structure determined by the condition \eqref{BandGenPos1} together with Eqs.~\eqref{f,app}-\eqref{w,app}; the positive spectrum has infinitely many gaps.
		\item \label{th3} $P_\sigma(H^{+R}_{ \ell_3,\ell_4 })=0$ for any $\ell_3,\ell_4>0\,$, and the spectrum is dominated by gaps; at high energies, the bands appear only in the vicinity of the numbers $k=\frac{n\pi}{a}$, $(n-\frac12)\frac{\pi}{a}$ and $\frac{(\sqrt{3}+1)n\pi}{2a}$, $n\in\mathbb{N}$.
		\item \label{th4} $P_\sigma(H^{+R}_ {0,\ell_4 })\approx 0.82\,$ for any $\ell_4>0$; in the limit $\ell_3\to 0$, the $R$ coupling at vertices of degree `three' tends to the Kirchhoff one resulting in wider spectral bands.
		\item \label{th5} The negative spectrum consists of at most three number of bands which may merge for particular values of parameters; there is no flat bands in this part of the spectrum.
		\item \label{th6} The negative bands become exponentially narrow as $a\to\infty$ shrinking to the eigenvalues of star graphs of degree three $(-3\ell_3^{-2})$ and four $(-\ell_4^{-2})$.
	\end{enumerate}
\end{theorem}

\section {Cairo model with the $(-1)^{d_v}R$ coupling}
\label{sect:RminusR}

Let us modify now the model by changing the matching condition in the vertices of degree three to $-R$. According to the discussion in Sec.~\ref{s:vertex}, there is no asymptotic decoupling in that case. As before, we will discuss different parts of the spectrum separately, and summarize the results in Theorem~\ref{thmGenR-R} at the end of the section.

\subsection{Positive spectrum}	
\label{sect:PosGenRminusR}

The secular equation is derived in the same way as above: we employ appropriate Ansätze for the wave function components and match them at the graph vertices, this time with the $-R$ coupling at vertices of degree three, considering again the clockwise orientation at all the vertices. The obtained spectral condition has the same structure as \eqref{spectGenPos}, however, with different functions $f(k),\, g(k),\, h(k)$ and $w(k)$ given explicitly in Appendix~\ref{sect:appB} through eqs. \eqref{f,minusR}--\eqref{w,minusR}. The resulting positive spectrum consists again of two parts:

\begin{itemize}
\setlength{\itemsep}{-3pt}
\item  If $\ell_3=\ell_4$, the number $k=\ell_3^{-1}$ belongs to the spectrum, this time for $a=\frac{(2n-1)\pi}{\sqrt{3}+3}\ell_3$ with $n\in\mathbb{N}$. To identify this flat band, one has to find when the functions in Eqs. \eqref{f,minusR}--\eqref{w,minusR} can vanish simultaneously. We start again with the simplest case, the function $h(k)$, which vanishes for three sets of parameters, specifically $k=1/\sqrt{\ell_3\,\ell_4}$ and $a=(n-\frac12)\pi\sqrt{\ell_3\,\ell_4}$, or $k=\frac{n\pi}{a}$ and $\ell_3=\ell_4$, or $k=\ell_3^{-1}$ and $\ell_3=\ell_4$. Inspecting then $g(k)$, $w(k)$ and $f(k)$ at these values, we see that the first two functions vanish at all the three sets, while the values of the third one, $f(k)$, are reduced respectively to $-1024 (-1)^n \ell_3^2 \ell_4^{-2} \big(4 \sin  \sqrt{3} (n-\frac{1}{2} )\pi    +(-1)^n  (\cos  \sqrt{3} (2 n-1)\pi -3 )\big)$ and  $1024 \pi ^8 n^8 \ell_3^8 a^{-8}\big(4 (-1)^n \cos  \sqrt{3} n\pi  +\cos  2 \sqrt{3} n\pi  +3\big)$ in the first two cases and $8192 \cos^4 \frac{ (\sqrt{3}+3 ) a}{2 \ell_3}$ is the last one. The former two expressions can be checked to be nonzero in view of the irrationality of $\sqrt{3}$, the last one verifies the claim.
\item  Apart from the flat bands mentioned above, the rest of the spectrum is absolutely continuous having a band-gap structure. Repeating the argument of Sec.~\ref{sect:PosGen}, the spectral bands are given by the conditions \eqref{BandGenPos1}, this time together with Eqs. \eqref{f,minusR}--\eqref{w,minusR}. The band-gap pattern in dependence on the parameter $a$ and the length scales $\ell_3$ and $\ell_4$ is illustrated in Figs.~\ref{Cairo,ell3=ell4=1b}--\ref{Cairo,ell3=1,a=1b} in which, for the sake of comparison, we use the same values as in the model with $R$ coupling at all vertices shown in  Figs.~\ref{Cairo,ell3=ell4=1a}--\ref{Cairo,ell3=1,a=1a}.
\end{itemize}

\subsection{High-energy regime and the probability $P_{\sigma}(H)$}
\label{sect:HighRminusR}

As we have mentioned, there is no high-energy decoupling of the edges in this case, the scattering matrix at vertices of degree four being given by \eqref{splus4}, while for vertices of degree three we now get from \eqref{sij,onshell} the limit
\begin{equation}\label{hightenS} 
 \lim_{k\to\infty} S^{- R}(k)=
 \frac13 \left(\begin{array}{ccc}
 	\phantom{-}1 &     -2 &       -2   \\
 	-2 &     \phantom{-}1 &       -2   \\
 	-2 &            -2 &        \phantom{-}1   \\
 \end{array}\right)  .
\end{equation}
As a result, the transport is non-trivial in this case in contrast to the model of the previous section and the probability of belonging to the spectrum, $P_{\sigma}(H)$, is expected to be nonzero. To determine it, we rewrite the spectral condition keeping the leading order of $k$ as $k\to\infty$ in the form
 \begin{align} \label{asymptcond,minusR}
 	4 \;\ell_4^4\;\ell_3^8\; \sin ^2 ka\,&\big\{ \tilde{f}(k)-\tilde{g}(k) (\cos \theta_1 +\cos \theta_2 )-\tilde{h}(k) (\cos 2 \theta_1 +\cos 2 \theta_2 )-\tilde{w}(k) \,\cos \theta_1   \,\cos \theta_2 \big\}\, k^{12} \\ \nonumber
 	&+\mathcal{O}(k^{10})=0,
 \end{align}
where
 \small
 \begin{equation}\label{fghwtildeminusR}
 	\begin{aligned}
 		\tilde{f}(k)&= 54 \cos 2 ka+60 \cos 4 ka+18 \cos 6 ka-\cos 2 (\sqrt{3}-4) ka-27 \cos 2 \sqrt{3} ka+60 \\
 		&  +16 \cos   (\sqrt{3}-3 ) ka -6 \cos  2  (\sqrt{3}-3 ) ka -11 \cos  2  (\sqrt{3}-2 ) ka +32 \cos(\sqrt{3}-1) ka \\
 		&  -12 \cos  2  (\sqrt{3}-1 ) ka -48 \cos   (\sqrt{3}+1 ) ka\ -54 \cos  2  (\sqrt{3}+1 ) ka -81 \cos 2 (\sqrt{3}+2 ) ka , \\[5pt]
 		\tilde{g}(k)&= 8 \cos 2 ka-8 \cos (\sqrt{3}-5) ka-32 \cos   (\sqrt{3}-3 ) ka -36 \cos  2 \sqrt{3} ka +56  \\
 		&  -4 \cos  2  (\sqrt{3}-2 ) ka -32 \cos   (\sqrt{3}-1 ) ka-24 \cos  2  (\sqrt{3}-1 ) ka +72 \cos   (\sqrt{3}+3 ) ka  ,\\[5pt]
 		\tilde{h}(k)&=-16 \sin ^2 ka, \\[5pt]
 		\tilde{w}(k)&=40 \cos 2 ka-16 \cos  (\sqrt{3}-3 ) ka-4 \cos  2  (\sqrt{3}-2 ) ka -36 \cos  2 \sqrt{3} ka +24 \\
 		& -32 \cos   (\sqrt{3}-1 ) ka -24 \cos  2  (\sqrt{3}-1 ) ka +48 \cos   (\sqrt{3}+1 ) ka  .
 	\end{aligned}
 \end{equation}
 \normalsize
Inspecting relation \eqref{asymptcond,minusR}, we see that there are narrow bands in the vicinity of the roots of $\sin ^2 ka$ which we are not going to specify, and wide bands corresponding to vanishing of the expression in the curly brackets, which give rise to the probability $P_{\sigma}(H)$. To find it, we repeat the argument of Sec.~\ref{sect:PosEll30High} regarding $x:=\sqrt{3}\, ka\,(\mathrm{mod}\,2\pi)$ and $y:=ka \,(\mathrm{mod}\,2\pi)$ as a pair of independent identically distributed random variables on $[0,2\pi)$. The mentioned bands are then determined again by the condition \eqref{con,ell30,xy1} in which the functions $\tilde{f}(x,y)$, $\tilde{g}(x,y)$, $\tilde{h}(x,y)$ and $\tilde{w}(x,y)$ are now given by
 		\small
 	\begin{equation}\label{fghtilde,xy,minusR}
 		\begin{aligned}
 			\tilde{f}(x,y) &=  2  (21 \sin 2y+35 \sin 4y-3 \sin 6y)\sin 2x- (66 \cos 2y+92 \cos 4y+6 \cos 6y+27)\cos 2x  \\
 			&    +2 \left(40 \sin x \sin y-8 \cos x \cos y+8 \cos (x-3 y)+12 (3 \cos 2y+5) \cos ^2 2y\right)-\cos 2 (x-4 y)                    , \\[5pt]
 			\tilde{g}(x,y) &=   -4 \left(\cos 2 (x-2 y)+8 \cos (x-3 y)+2 \cos (x-5 y)+9 \cos 2 x\right)         \\
 			&  +8 \left(-4 \cos (x-y)-3 \cos 2 (x-y)+9 \cos (x+3 y)+\cos 2 y+7\right)                     , \\[5pt]
 			\tilde{h}(x,y) &=-16 \sin ^2 y, \\[5pt]
 			\tilde{w}(x,y) &=-4 \left(6 \cos 2 (x-y)+\cos 2 (x-2 y)+4 \cos (x-3 y)+9 \cos 2 x-6\right)  \\
 			& -4 \left(20 \sin x \sin y-4 \cos x \cos y-10 \cos 2 y\right).
 		\end{aligned}
 	\end{equation}
\normalsize
The probability \eqref{probsigma} can be determined numerically by finding the area of the region in which the condition \eqref{con,ell30,xy1} with the new functions \eqref{fghtilde,xy,minusR} holds on $[0,2\pi)\times[0,2\pi)$. It is shown in Fig.~\ref{FigProbRminusR}; the resulting value is $P_{\sigma}(H)\approx 0.82$.

We note that it coincides with the probability obtained for the Kirchhoff limit in the previous section. This conclusion is not only numerical. To see that, we observe that the functions in \eqref{fghtilde,xy} and \eqref{fghtilde,xy,minusR} are obtained one from the other through the transformations $x\leftrightarrow x+\frac{\pi}2$ and $y\leftrightarrow y-\frac{\pi}2$, in other words, that that the two regions in which the band condition \eqref{con,ell30,xy1} together with any of the sets of functions \eqref{fghtilde,xy} and \eqref{fghtilde,xy,minusR} holds coincide on the squares $[0,2\pi)\times[0,2\pi)$ and $[\frac{\pi}2,5\frac{\pi}2)\times [-\frac{\pi}2,3\frac{\pi}2)$ -- cf. Fig~\ref{FigProbsArea} -- which together with the $2\pi$-periodicity in both the $x$ and $y$ directions yields the claim. This exact coincidence is not surprising in view of the fact that the high-energy limit of the scattering matrix \eqref{hightenS} is the same as obtained from the Kirchhoff coupling.

\subsection{Negative spectrum}	
\label{sect:NegGenRminusR}

We have to replace $k$ by $i\kappa$ with $\kappa>0$ again, the spectral condition is then given by \eqref{BandGenNeg1} in which $f(\kappa)$ as well as $g(\kappa),\, h(\kappa)$ and $w(\kappa)$ that are included in $\Lambda_i(\kappa)$ are now given by Eqs.~\eqref{f,minusR}--\eqref{w,minusR} with the indicated replacement. The negative spectrum, consisting of at most three bands, is illustrated in Figs.~\ref{Cairo,ell3=ell4=1b}--\ref{Cairo,ell3=1,a=1b}. We see that:

\begin{itemize}
\setlength{\itemsep}{-3pt}
\item For large values of the parameter $a\,$, the negative bands shrink to points, cf. Figs.~\ref{Cairo,ell3=ell4=1b} and \ref{Cairo,ell3=2,ell4=23b}, corresponding to the eigenvalues of star graphs of degree four and three, $E_{+R}=-\ell_4^{-2}$ and $E_{-R}=-\tfrac13\ell_3^{-2}$, respectively -- cf.~\eqref{Neg,evRStarG}. This can be easily seen by rewriting the spectral condition in the form
\begin{equation}\label{Neg,gen,leadminus}
 \left(1-\kappa ^2 \ell_4^2\right)^2 \left(1-3\kappa ^2 \ell_3^2\right)^4 \, \e^  {2 (\sqrt{3}+3) \kappa a}+\mathcal{O}(\e^  {8 \kappa a})=0.
\end{equation}
To find the asymptotic band widths, we set $\kappa_1=\ell_4^{-1}+\delta$ and $\kappa_2=(\sqrt{3}\,\ell_3)^{-1}+\delta\,$; in analogy with Sec.~\ref{sect:NegGen} we solve the resulting linear/quadratic equation for $\delta$ in the two cases, respectively, arriving at
$$ \triangle E_1=\frac{8 \ell_3^2}{\ell_4^4-3 \ell_4^2 \ell_3^2}\,\e^{\frac{-2 a}{\ell_4}}+\mathcal{O}\big(\e^{\frac{- (\sqrt{3}+1) a}{\ell_4}}\big),$$
and
$$\triangle E_2=\frac{8 \sqrt{2} \left(\ell_4^2+3 \ell_3^2\right)}{9 \ell_3^2 \left(3 \ell_3^2-\ell_4^2\right)}\,\e^{\frac{-2 a}{\sqrt{3}\,\ell_3}}+\mathcal{O}\big(\e^{\frac{- (\sqrt{3}+3) a}{3\,\ell_3}}\big).$$
\item There is no flat band; as in Sec.~\ref{sect:NegGen} one has to check that the coefficients of the quasimomentum-dependent terms, the functions $g(\kappa)$, $h(\kappa)$ and $w(\kappa)$, cannot vanish simultaneously. Focusing on the simplest of them, the condition $h(\kappa)=0$ can be simplified to
$$\rho(\kappa):=(\kappa ^2 \ell_4^2+1) (\kappa ^2 \ell_3^2+1) \cosh 2 \kappa a-\left(\kappa ^4 \ell_4^2 \,\ell_3^2-\kappa ^2 (\ell_4^2-4 \ell_4 \,\ell_3+\ell_3^2)+1\right)=0,$$
and since $\cosh x>1$ for $x>0$, we have $\rho(\kappa)> 2 \kappa ^2 (\ell_4-\ell_3)^2$ and since the inequality is sharp, the claim is verified. 
\end{itemize}

\subsection{The Dirichlet limit, $\ell_3\to 0$}
\label{sect:Ell30minusR}

Let us consider again the limit $\ell_3\to 0$, now for the $-R$ coupling at vertices of degree three. Summing the condition \eqref{couplings} in this case over $j$, we get $\sum_j \psi_j=0$, and in the limit $\ell_3\to 0$ the condition becomes $\psi_{j+1}=-\psi_j$ which in the case of an odd vertex degree means $\psi_j=0$. Consequently, the Cairo lattice decomposes into an infinite family of finite graphs (consisting of one or four edges) with Dirichlet boundary the spectrum of which is pure point and infinitely degenerate.

\subsubsection{Positive spectrum}	
\label{sect:PosEll3030minusR}

All the three functions \eqref{g,minusR}--\eqref{w,minusR} in \eqref{spectGenPos} have a multiplicative factor $\ell_3$, hence the quasimomentum-dependent part vanishes in the limit and the spectrum consists of flat bands only. More explicitly, the condition $f(k)=0$ simplifies to

\begin{equation}\label{posell30,minusR}
	 32  \left[(k^2 \ell_4^2+1) \cos 2 ka+k^2 \ell_4^2-1\right]^2 \,\sin ^2 2 ka \;\sin ^2 (\sqrt{3}-1) ka=0 ,
\end{equation}
so the eigenvalues come from zeros of the three factors, being the squares of $k=\frac{n\pi}{(\sqrt{3}-1)a}$ and $\frac{n\pi}{2a}$ with $n\in\mathbb{N}$, and of the roots of the expression in the square brackets. More precisely, the latter are given by the condition $\frac{1-k^2 \ell_4^2}{1+k^2 \ell_4^2}=\cos 2 ka$ in which the left-hand side is decreasing with respect to $k$ on the interval $k\in(0,\infty)$, ranging from $1$ to $-1$, giving rise to an infinite number of eigenvalues. The three cases can be easily identified, the second one are the Dirichlet eigenvalues on an interval of length $b$, the other two refer to the four-legged cross of edge length $a$, the first belonging to Dirichlet eigenfunctions vanishing in the center of the cross and last one to the `genuine' eigenfunctions of the cross with the $R$ coupling in the center. Needles to say, $P_{\sigma}(H)=0$ holds in this case.

\subsubsection{Negative spectrum}	
\label{sect:NegEll3030minusR}
Concerning the negative eigenvalues, replacing $k$ by $i\kappa$ in \eqref{posell30,minusR}, we arrive at the spectral condition
\begin{equation}\label{negell30,minusR}
32  \left[(1-\kappa ^2 \ell_4^2) \cosh 2 \kappa a -\kappa ^2 \ell_4^2-1\right]^2 \, \sinh ^2  2 \kappa a  \; \sinh ^2(\sqrt{3}-1) \kappa a =0\,;
\end{equation}
since $\sinh x>0$ for $x>0$, it remains to check the term in the square brackets; clearly, a `genuine' cross eigenfunction only may contribute to the negative spectrum. We claim that an isolated negative eigenvalue may correspond to a $\kappa\in(0,\ell_4^{-1})$, and that happens if and only if $a>\ell_4$.

\smallskip

Indeed, let us denote the square-bracket expression as $\mathcal{G}(\kappa)$, obviously nonzero at $\kappa=\ell_4^{-1}$. Introducing positive $\omega:=\kappa\, \ell_4$ and $\lambda:=a\ell_4^{-1}$, the condition $\mathcal{G}(\kappa)=0$ reads $\frac{1+\omega ^2 }{1-\omega ^2}=\cosh 2 \lambda \omega $. Since the right-hand side is positive, the equation may have solutions only for $\omega<1$, or equivalently, $0<\kappa<\ell_4^{-1}$, for which the left-hand side of the equation is always positive. On the other hand, the equation has a unique solution for $\lambda$ which is $\lambda=\frac1{2\omega}\,\arccosh\frac{1+\omega ^2 }{1-\omega ^2}>1$. The last claim, $\lambda>1$, can be proved by showing that $\xi(\omega):=\arccosh\frac{1+\omega ^2 }{1-\omega ^2}-2\omega$ is positive. Computing the first derivative, we get $\xi'(\omega)= \frac{2\omega ^2}{1-\omega ^2}$ which is positive for $\omega<1$, we see that $\xi(\omega)$ is monotonically increasing on the domain; in combination with $\xi(0)=0$ and $\lim_{\omega\to 1} \xi(\omega)=+\infty$ this concludes the argument.

\subsection{The results summary}
\label{sect:mainRminusR}

Let us finally summarize the obtained results: let $H^{\pm R}_{ \ell_3,\ell_4 }$ be the quantum graph Hamiltonian described in the beginning of this section, then our observations can be summarized as follows:
\begin{theorem}\label{thmGenR-R}
	The spectrum of $H^{\pm R}_{ \ell_3,\ell_4 }$ has the following properties:
	\begin{enumerate}[label=\textnormal{(\roman*)}]
		\setlength{\itemsep}{-3pt}
		\item \label{th1R-R} The number $k=\ell_3^{-1}$ for $\ell_3>0$ always belongs to $\sigma(H^{\pm R}_{ \ell_3,\ell_3  })$ when $a=\frac{(2 n-1)\pi  }{\sqrt{3}+3}\ell_3$ with $n\in\mathbb{N}$.
		\item \label{th2R-R} Apart from the above infinitely degenerate eigenvalues, the spectrum of $H^{\pm R}_{ \ell_3,\ell_4 }$ for $\ell_3,\ell_4>0\,$ is absolutely continuous having a band-and-gap structure determined by the condition \eqref{BandGenPos1} together with Eqs. \eqref{f,minusR}--\eqref{w,minusR}; the positive spectrum has infinitely many gaps.
		\item \label{th3R-R} $P_\sigma(H^{\pm R}_{ \ell_3,\ell_4 })\approx 0.82\,$ holds for any $\ell_3,\ell_4>0\,$; vertices of both the even and odd parity remain transparent at high energies.
		\item \label{th4R-R} In the limit $\ell_3\to 0$ we have $P_\sigma(H^{\pm R}_{ 0,\ell_4  })=0\,$ for any $\ell_4>0$; the $-R$ coupling at the vertices of odd degree produces Dirichlet decoupled edges resulting in a pure point spectrum, infinitely degenerate, consisting of the squared roots of the functions $\sin 2 ka \,$, $\sin (\sqrt{3}-1) ka$ and $(k^2 \ell_4^2+1) \cos 2 ka+k^2 \ell_4^2-1 $.
		\item \label{th5R-R} The negative spectrum of $H^{\pm R}_{ \ell_3,\ell_4 }$ with $\ell_3,\ell_4>0\,$ consists of at most three number of bands which may merge for particular values of parameters; there is no negative flat band. The negative spectrum of $H^{\pm R}_{ 0,\ell_4  }$, $\ell_4>0\,$ consists of at most one isolated eigenvalue referring to $\kappa\in(0,\ell_4^{-1})$; it exists if and only if $a>\ell_4$.
       \item \label{th6-R} The negative bands of $H^{\pm R}_{ \ell_3,\ell_4 }$ with $\ell_3,\ell_4>0\,$ become exponentially narrow as $a\to\infty$ shrinking to the eigenvalues of star graphs of degree three $(-\tfrac13\ell_3^{-2})$ and four $(-\ell_4^{-2})$.
    \end{enumerate}
\end{theorem}


\appendix
\section*{Appendices}

\section{Explicit form of some functions appearing in Sec.~\ref{sect:PosGenRR}}
\label{sect:appA}
\renewcommand{\theequation}{\thesection.\arabic{equation}}
\setcounter{equation}{0}
The functions $f(k),\; g(k),\; h(k)$ and $w(k)$ entering the spectral condition \eqref{spectGenPos} in Sec.~\ref{sect:PosGen}, corresponding to the model with $R$ coupling at all the vertices, are as follows
\footnotesize
\begin{align}	\label{f,app}
	f(k)&=-2 \left(752 k^4 \ell_4^2 \ell_3^2+128 k^2 \ell_4 \ell_3 (k^2 \ell_4^2+1)+33 (k^2 \ell_4^2+1)^2+128 k^4 \ell_4 \ell_3^3 (k^2 \ell_4^2+1)-8 k^6 \ell_3^6 (k^4 \ell_4^4+1)\right) \\ \nonumber
	& +16 (k^2 \ell_3^2+3)  \left(k^6  (3 \ell_4^4 \ell_3^2+4 \ell_4^3 \ell_3^3)+k^4  (\ell_4^4+12 \ell_4^3 \ell_3+18 \ell_4^2 \ell_3^2+4\ell_4 \ell_3^3)+3 k^2 (2 \ell_4^2+4 \ell_4 \ell_3+\ell_3^2)+1\right)\cos(\sqrt{3}+3) ka\\ \nonumber
	& -16 (k^2 \ell_3^2-1)^2  \left(k^4 \ell_4^3 (\ell_4+4 \ell_3)+2 k^2 \ell_4 (3 \ell_4+2 \ell_3)+1\right)\cos (\sqrt{3}-5)ka-(k^2 \ell_4^2+1)^2 (k^2 \ell_3^2-1)^4 \cos 2 (\sqrt{3}-5) ka \\ \nonumber
	& -(k^2 \ell_4^2+1)^2 (k^2 \ell_3^2+3)^4 \cos 2 (\sqrt{3}+3) ka+4 (k^4 \ell_4^4-1) (k^2 \ell_3^2-1)^4 \cos 2 (\sqrt{3}-4)ka-4 k^4 \ell_3^4 (17 k^4 \ell_4^4+42 k^2 \ell_4^2+17) \\ \nonumber
	& -4 (k^4 \ell_4^4-1) (k^2 \ell_3^2-1)^4 \cos 2 (\sqrt{3}-2) ka+96 (k^2 \ell_4^2-1)^2 (k^2 \ell_3^2-1)^2 \cos (\sqrt{3}-1)ka -2 k^8 \ell_3^8 (5 k^4 \ell_4^4-6 k^2 \ell_4^2+5) \\ \nonumber
	& +2 (k^2 \ell_3^2-1)^2  \left(k^2 \big(-7 k^2 \ell_4^4+k^2 \ell_3^4 (5 k^4 \ell_4^4-6 k^2 \ell_4^2+5)+2 \ell_3^2 (k^4 \ell_4^4-6 k^2 \ell_4^2+1)+82 \ell_4^2\big)-7\right)\cos 2 (\sqrt{3}-1)ka \\ \nonumber
	& +8 (k^2 \ell_4^2-1)^2 \left(k^8 \ell_3^8+k^6 \ell_3^6-3 k^4 \ell_3^4+7 k^2 \ell_3^2-6\right) \cos 4ka  -8 (k^4 \ell_4^4-1) (k^2 \ell_3^2-1)^2 (k^4 \ell_3^4+4 k^2 \ell_3^2+3) \cos 6ka \\ \nonumber
	& +8 (k^2 \ell_3^2-1)  \left(k^{10} \ell_4^4 \ell_3^6-k^8 \ell_4^4 \ell_3^4+k^6 (7 \ell_4^4 \ell_3^2-\ell_3^6)+k^4 (9 \ell_4^4+32 \ell_4^3 \ell_3+\ell_3^4)-k^2 \ell_3 (32 \ell_4+7 \ell_3)-9\right)\cos2ka \\ \nonumber
	& -4 (k^2 \ell_3^2-1) (k^2 \ell_3^2+3)^2 \left(k^2 \ell_3^2 (k^2 \ell_4^2-1)^2 \cos 2 (\sqrt{3}+1) ka-(k^4 \ell_4^4-1) (k^2 \ell_3^2+3) \cos 2 (\sqrt{3}+2)ka\right) \\ \nonumber
	& -4 k^2 \ell_3 (k^2 \ell_4^2-1) (k^2 \ell_3^2-1)^2 \left(\ell_3 (k^2 \ell_4^2-1) (k^2 \ell_3^2-1) \cos 2 (\sqrt{3}-3) ka-32 \ell_4 \cos (\sqrt{3}-3) ka\right) \\ \nonumber
	& -4 (k^4 \ell_4^4-1) (k^2 \ell_3^2-1)^3 (k^2 \ell_3^2+3) \cos 2 \sqrt{3} ka+2 (k^2 \ell_4^2+1)^2 (k^4 \ell_3^4+2 k^2 \ell_3^2-3)^2 \cos 8ka \\ \nonumber
	& -128 (k^2 \ell_4^2-1) (k^2 \ell_3^2-1)  \left(k^4 \ell_4 \ell_3^2 (\ell_4+\ell_3)+k^2 (\ell_4^2+3 \ell_4 \ell_3+\ell_3^2)+1\right)\cos (\sqrt{3}+1) ka. \nonumber &&
\end{align}


\begin{align}	\label{g,app}
g(k)&=8 (k^4 \ell_4^4-1) (k^2 \ell_3^2-1)^3 \cos (\sqrt{3}-7)ka-16 (k^4 \ell_4^4-1) (k^2 \ell_3^2-1)^2 (k^2 \ell_3^2+3) \cos 2 \sqrt{3} ka-104  \\ \nonumber
	& +8 (k^4 \ell_4^4-1) (k^2 \ell_3^2-1)^2 (5 k^2 \ell_3^2+3) \cos (\sqrt{3}-3) ka -8 (k^4 \ell_4^4-1) (5 k^6 \ell_3^6+k^4 \ell_3^4-k^2 \ell_3^2-5) \cos (\sqrt{3}+1) ka \\ \nonumber
	& +4(k^2 \ell_4^2-1)^2 (k^2 \ell_3^2-1) (k^2 \ell_3^2+3)^2 \cos 2 (\sqrt{3}+1) ka-8 (k^4 \ell_4^4-1) (k^2 \ell_3^2-1) (k^2 \ell_3^2+3)^2 \cos (\sqrt{3}+5) ka\\ \nonumber
	& +8 k^2 \left(\ell_4^2 (-2 k^6 \ell_3^6+6 k^4 \ell_3^4+58 k^2 \ell_3^2+2)-\ell_3^2 (3 k^4 \ell_3^4+7 k^2 \ell_3^2+9)-k^2 \ell_4^4  (3 k^6 \ell_3^6+7 k^4 \ell_3^4+9 k^2 \ell_3^2+13 )\right)\\ \nonumber
	& +32 (k^4 \ell_4^4-1) (k^6 \ell_3^6+3 k^4 \ell_3^4-k^2 \ell_3^2-3) \cos 2ka-8 (k^2 \ell_4^2-1)^2 (k^6 \ell_3^6+5 k^4 \ell_3^4-5 k^2 \ell_3^2-1) \cos 4ka \\ \nonumber
	& +16 (k^2 \ell_3^2+3)  \left(2 k^4 \ell_3^4 (k^4 \ell_4^4+1)+3 k^4 \ell_4^4-2 k^2 \ell_4^2-k^2 \ell_3^2 (k^4 \ell_4^4+6 k^2 \ell_4^2+1)+3\right)\cos (\sqrt{3}+3) ka \\ \nonumber
	&  -16 (k^2 \ell_3^2-1)^2  \left((2 k^2 \ell_3^2 (k^4 \ell_4^4+1)-(k^2 \ell_4^2+1)^2)\cos (\sqrt{3}-5)ka+2 (k^2 \ell_4^2-1)^2  \cos (\sqrt{3}-1) ka\right) \\ \nonumber
	& +4(k^2 \ell_4^2-1)^2 (k^2 \ell_3^2-1)^3 \cos 2 (\sqrt{3}-3) ka-16 (k^4 \ell_4^4-1) (k^2 \ell_3^2-1)^3 \cos 2 (\sqrt{3}-2) ka \\ \nonumber
    & +8 (k^2 \ell_3^2-1)^2  \left(k^2 (\ell_4^2 (k^2 (\ell_3^2 (3 k^2 \ell_4^2+2)+\ell_4^2)-10)+3 \ell_3^2)+1\right)\cos 2 (\sqrt{3}-1) ka . \nonumber &&
\end{align}


\begin{flalign}\label{h,app}
	 h(k)&=64 \left( (k^2 \ell_4 \ell_3-1)^2\,\sin ^2 ka+k^2 (\ell_4-\ell_3)^2 \cos ^2 ka\right)^2 .  &&
\end{flalign}


\begin{flalign}	\label{w,app}
	w(k)	&= -16  (k^2 \ell_3^2-1 )^2  \left(k^2 (\ell_4^2 (k^2 (\ell_3^2 (3 k^2 \ell_4^2+2)+\ell_4^2)-10)+3 \ell_3^2)+1\right)\sin ^2  (\sqrt{3}-1 )ka \\ \nonumber
	&  -16 (k^2 \ell_3^2+3)  \left(k^6 \ell_4^3 \ell_3^2 (3 \ell_4-4 \ell_3)+3 k^2 (2 \ell_4^2-4 \ell_4 \ell_3+\ell_3^2)+k^4 \ell_4 (\ell_4^3-12 \ell_4^2 \ell_3+18 \ell_4 \ell_3^2-4 \ell_3^3)+1\right) \cos (\sqrt{3}+3)ka  \\ \nonumber
	&  +4(k^2 \ell_4^2-1)^2 (k^2 \ell_3^2-1) (k^2 \ell_3^2+3)^2 \cos 2 (\sqrt{3}+1) ka   +128 k^2 \ell_4 \ell_3 (k^2 \ell_4^2-1) (k^2 \ell_3^2-1)^2 \cos (\sqrt{3}-3) ka \\ \nonumber
	&  -16 (k^4 \ell_4^4-1) (k^2 \ell_3^2-1)^2 (k^2 \ell_3^2+3) \cos 2 \sqrt{3}ka-8 (k^2 \ell_4^2-1)^2 (k^6 \ell_3^6+k^4 \ell_3^4+3 k^2 \ell_3^2-5) \cos 4ka   \\ \nonumber
	&  +128 (k^2 \ell_4^2-1) (k^2 \ell_3^2-1)  \left(k^4 \ell_4 \ell_3^2 (\ell_4-\ell_3)+k^2 (\ell_4^2-3 \ell_4 \ell_3+\ell_3^2)+1\right) \cos (\sqrt{3}+1) ka  \\ \nonumber
	&  +4(k^2 \ell_4^2-1)^2 (k^2 \ell_3^2-1)^3 \cos 2 (\sqrt{3}-3)ka+32 (k^4 \ell_4^4-1) (k^4 \ell_3^4-1) (k^2 \ell_3^2-1) \cos 2ka   \\ \nonumber
	&  -96 (k^2 \ell_4^2-1)^2 (k^2 \ell_3^2-1)^2 \cos (\sqrt{3}-1)ka-16 (k^4 \ell_4^4-1) (k^2 \ell_3^2-1)^3 \cos 2 (\sqrt{3}-2) ka   \\ \nonumber
	&  +4(k^2 \ell_3^2-1)^2  \left(4 k^2 \ell_4 (\ell_4 (k^2 \ell_4 (\ell_4-4 \ell_3)+6)-4 \ell_3)+4\right)  \cos (\sqrt{3}-5)ka  . \nonumber &&
\end{flalign}


\noindent Furthermore, the functions $\mathcal{F}_j,\: j=1,2,3,$ appearing in Sec.~\ref{sect:PosGenHigh}(i) are of the form

\begin{align*}
	\mathcal{F}_1	&=   \sqrt{\frac{-6 \sqrt[3]{2} \,\zeta_1\, \zeta_3+2 \sqrt[3]{2} \,\zeta_2^2-2 \,\zeta_2 \sqrt[3]{\zeta_4}+ (2\,\zeta_4)^{2/3}}{6 \,\zeta_1\, \sqrt[3]{\zeta_4}}},    \\[6pt]
\mathcal{F}_2	&=   \sqrt{\frac{\left(6 \sqrt[3]{2}  \,\zeta_1\, \zeta_3-2 \sqrt[3]{2} \, \zeta_2^2\right)\left(1-i\sqrt{3}\right) -4 \zeta_2 \,\sqrt[3]{\zeta_4}- \left(1+i \sqrt{3}\right) (2\,\zeta_4)^{2/3}}{12 \,\zeta_1\, \sqrt[3]{\zeta_4}}}   ,   \\[6pt]
\mathcal{F}_3	&=   \sqrt{\frac{\left(6 \sqrt[3]{2}  \,\zeta_1\, \zeta_3-2 \sqrt[3]{2}   \,\zeta_2^2\right)\left(1+i\sqrt{3}\right)-4 \zeta_2\, \sqrt[3]{\zeta_4}- \left(1-i\sqrt{3}\right) (2\,\zeta_4)^{2/3}}{12 \,\zeta_1\, \sqrt[3]{\zeta_4}}}   ,
\end{align*}
where
\begin{align*}
  \zeta_1	&=  \frac{1}{a^6}\left(   512 \, \pi ^{12} \, \ell_4^4\, \ell_3^8 \, \sin ^2  \sqrt{3} n\pi     \right)  ,     \\
\zeta_2	&=  - \frac{1}{a^6}\left(   128\,  \pi ^{10}\,  \ell_4^2\,  \ell_3^6  \left(\ell_4^2 (\tau_1+9)+8 \ell_3^2\right)  \sin ^2 \sqrt{3} n\pi  \right) ,    \\
\zeta_3	&=    \frac{1}{a^6}\left(   256 \, \pi ^8 \, \ell_3^4 \left(\ell_4^4 (\tau_1+3)-\ell_4^2\,  \ell_3^2 (\tau_1-3)+2 \ell_3^4\right)  \sin ^2 \sqrt{3} n\pi  \right),    \\
\zeta_4	&=  3 \sqrt{3\,  \zeta_1^2 \left(27 \, \zeta_1^2\,  \zeta_5^2-18\,  \zeta_1 \, \zeta_2\,  \zeta_3 \, \zeta_5+4 \zeta_1 \, \zeta_3^3+4 \zeta_2^3\,  \zeta_5-\zeta_2^2 \, \zeta_3^2\right)}-27 \zeta_1^2 \, \zeta_5+9 \zeta_1 \, \zeta_2 \, \zeta_3-2 \zeta_2^3 , \\
\zeta_5	&=   \frac{1}{a^6}\left(   64 \pi ^6 (\tau_1+1) \left(\ell_3^3-\ell_4^2\,  \ell_3\right)^2 \left(\cos 2 \sqrt{3} n\pi\ -1\right)   \right).
\end{align*}

\normalsize
\section{Explicit form of some functions appearing in Sec.~\ref{sect:RminusR}}
\label{sect:appB}
\renewcommand{\theequation}{\thesection.\arabic{equation}}
\setcounter{equation}{0}
The functions $f(k),\; g(k),\; h(k)$ and $w(k)$ entering the spectral condition \eqref{spectGenPos} in Sec.~~\ref{sect:RminusR} are as follows
\footnotesize
\begin{flalign}	\label{f,minusR}
	f(k)	&=  2 \left(5 k^4 \ell_4^4-6 k^2 \ell_4^2-8 k^2 \ell_3^2 (k^4 \ell_4^4+1)+2 k^4 \ell_3^4 (17 k^4 \ell_4^4+42 k^2 \ell_4^2+17)+128 k^6 \ell_4 \ell_3^5 (k^2 \ell_4^2+1)+5\right) \\ \nonumber
	&  +2 (k^2 \ell_3^2-1)^2 \left(k^2 (-5 k^2 \ell_4^4+k^2 \ell_3^4 (7 k^4 \ell_4^4-82 k^2 \ell_4^2+7)-2 \ell_3^2 (k^4 \ell_4^4-6 k^2 \ell_4^2+1)+6 \ell_4^2)-5\right)\cos  2 (\sqrt{3}-1) ka   \\ \nonumber
	&    +(k^2 \ell_3^2-1)^4 \left((k^2 \ell_4^2+1)^2 \cos 2 (\sqrt{3}-5) ka+8 (k^4 \ell_4^4-1)   \sin 2ka \,\sin 2 (\sqrt{3}-3)ka\right)+256 k^8 \ell_4 \ell_3^7(k^2 \ell_4^2+1)\\ \nonumber
	&  +32k^4\ell_3^3 (k^2 \ell_4^2-1)(k^2 \ell_3^2-1)^2\left(3  \ell_3 (k^2 \ell_4^2-1)  \cos  (\sqrt{3}-1) ka-4 \ell_4   \cos (\sqrt{3}-3) ka\right) +66 k^8 \ell_3^8 (k^2 \ell_4^2+1)^2\\ \nonumber
	&   -8 (k^2 \ell_3^2-1) \left(9 k^{10} \ell_4^4 \ell_3^6+k^8 \ell_4^3 \ell_3^4 (7 \ell_4+32 s)-k^6 \ell_3^2 (\ell_4^4+32 \ell_4 \ell_3^3+9 \ell_3^4)+k^4 (\ell_4^4-7 \ell_3^4)+k^2 \ell_3^2-1\right) \cos 2 ka   \\ \nonumber
	&   -128 k^4 \ell_3^3 (k^2 \ell_4^2-1) (k^2 \ell_3^2-1)  \left(3 k^2 \ell_4 \ell_3^2+k^2 \ell_3^3+\ell_4^2 (k^4 \ell_3^3+k^2 \ell_3)+\ell_4+\ell_3\right)\cos (\sqrt{3}+1) ka+1504 k^8 \ell_4^2 \ell_3^6   \\ \nonumber
	&   +8 (k^2 \ell_4^2-1)^2 \left(6 k^8 \ell_3^8-7 k^6 \ell_3^6+3 k^4 \ell_3^4-k^2 \ell_3^2-1\right) \cos 4 ka-2 (k^2 \ell_4^2+1)^2 (-3 k^4 \ell_3^4+2 k^2 \ell_3^2+1)^2 \cos 8 ka   \\ \nonumber
	&    +4 (k^4 \ell_4^4-1) (k^2 \ell_3^2-1)^3 (3 k^2 \ell_3^2+1) \cos 2 \sqrt{3} ka-4 (k^4 \ell_4^4-1) (k^2 \ell_3^2-1) (3 k^2 \ell_3^2+1)^3 \cos 2 (\sqrt{3}+2) ka \\ \nonumber
	&  -4 (k^2 \ell_4^2-1)^2 (k^2 \ell_3^2-1) (3 k^2 \ell_3^2+1)^2 \cos 2 (\sqrt{3}+1) ka +(k^2 \ell_4^2+1)^2 (3 k^2 \ell_3^2+1)^4 \cos2 (\sqrt{3}+3) ka \\ \nonumber
	&  -4 (k^2 \ell_4^2-1)^2  (k^2 \ell_3^2-1 )^3 \cos 2 (\sqrt{3}-3) ka-8 (k^4 \ell_4^4-1) (k^2 \ell_3^2-1)^2 \left(3 k^4 \ell_3^4+4 k^2 \ell_3^2+1\right) \cos 6 ka  \\ \nonumber
	& +16 k^4 \ell_3^3 (3 k^2 \ell_3^2+1)  \left(	
	(12 k^2  \ell_3^2+4) (k^2 \ell_4^3+\ell_4)+(k^2 \ell_3^3+ 3 \ell_3)(k^4 \ell_4^4+6 k^2 \ell_4^2+1)
	\right) \cos (\sqrt{3}+3) ka \\ \nonumber
	&  -16 k^4 \ell_3^3 (k^2 \ell_3^2-1)^2  \left(\ell_4 (k^2 \ell_4 (\ell_4 (k^2 \ell_4\ell_3 +4)+6 s)+4)+s\right)\cos (\sqrt{3}-5) ka . \nonumber &&
\end{flalign}

%

\begin{flalign}	\label{g,minusR}
	g(k)	&=4 k^2 \ell_3^2\times\Bigg[      2 k^2 \left(k^2 (-\ell_4^4  (k^2 \ell_3^2-1 )  (3 k^2 \ell_3^2+1 )^2+9 k^2 \ell_3^6-3 \ell_3^4)-5 \ell_3^2\right)  \cos  (\sqrt{3}+5) ka \\ \nonumber
	&  +4 k^2  \left(k^2 (\ell_4^4 (k^2 \ell_3^2-1)^3-k^2 \ell_3^6+3 \ell_3^4)-3 \ell_3^2\right)\cos  2 (\sqrt{3}-2) ka-2 (k^2 \ell_4^2-1)^2 (k^6 \ell_3^6+5 k^4 \ell_3^4-5 k^2 \ell_3^2-1) \cos 4 ka  \\ \nonumber
	&   -4 k^2  \left(k^2 (\ell_4^4 (k^2 \ell_3^2-1)^2 (3 k^2 \ell_3^2+1)-3 k^2 \ell_3^6+5 \ell_3^4)-\ell_3^2\right)\cos 2 \sqrt{3} ka -8 k^2 \ell_3^2  (k^2 \ell_4^2-1 )^2  (k^2 \ell_3^2-1 )^2 \cos   (\sqrt{3}-1 ) ka\\ \nonumber
	&   +2 k^2  \left(k^2 (\ell_4^4 (k^2 \ell_3^2-1)^3-k^2 \ell_3^6+3 \ell_3^4)-3 \ell_3^2\right)  \cos  (\sqrt{3}-7 ) ka  -8 (k^4 \ell_4^4-1) (3 k^6 \ell_3^6+k^4 \ell_3^4-3 k^2 \ell_3^2-1) \cos 2 ka  \\ \nonumber
	&    +2 k^2 \left(\ell_3^2 (13 k^4 \ell_3^4+9 k^2 \ell_3^2+7)+k^2 \ell_4^4 (13 k^6 \ell_3^6+9 k^4 \ell_3^4+7 k^2 \ell_3^2+3)-2 \ell_4^2 (k^6 \ell_3^6+29 k^4 \ell_3^4+3 k^2 \ell_3^2-1)\right)+6   \\ \nonumber
	& + 2\left( -8 \sin ^4 ka- k^6 \ell_3^6- k^4 \ell_3^4- k^2 \ell_4^2  (k^2 \ell_3^2-1 )^2  (k^2  (\ell_4^2  (k^2 \ell_3^2+3 )-10 \ell_3^2 )+2 )+5 k^2 \ell_3^2  \right) \cos 2 (\sqrt{3}-1) ka  \\ \nonumber
	&   +k^2 \left(-k^2 \ell_4^4-3 k^2 \ell_3^4  (k^2 \ell_4^2-1 )^2-5 \ell_3^2  (k^2 \ell_4^2-1 )^2+9 k^4 \ell_3^6  (k^2 \ell_4^2-1 )^2+2 \ell_4^2\right)  \cos 2 (\sqrt{3}+1) ka\\ \nonumber
	&  +4 k^2  \left( k^2  (9 k^2 \ell_3^6+\ell_4^4  (9 k^6 \ell_3^6+5 k^2 \ell_3^2+2 )-2 \ell_4^2 \ell_3^2  (3 k^4 \ell_3^4+10 k^2 \ell_3^2+3 ) )+5 \ell_3^2 \right) \cos (\sqrt{3}+3) ka\\ \nonumber
	&   +4 k^2 \left(k^2  (\ell_4^4 (k^2 \ell_3^2-2)  (k^2 \ell_3^2-1 )^2+2 \ell_4^2 \ell_3^2  (k^2 \ell_3^2-1 )^2+k^2 \ell_3^6-4 \ell_3^4 )+5 \ell_3^2\right) \cos  (\sqrt{3}-5) ka  \\ \nonumber
	& -2\left(  16 \sin ^4 ka+ k^2 \ell_3^2+ k^6  (\ell_4^4 (5 k^6 \ell_3^6+k^4 \ell_3^4-k^2 \ell_3^2-5)-\ell_3^4 (5 k^2 \ell_3^2+1) )  \right) \cos(\sqrt{3}+1) ka \\ \nonumber
	&   + 2\left(  16 \sin ^4 ka-7 k^2 \ell_3^2- k^4  (\ell_4^4 (k^2 \ell_3^2-1)^2 (3 k^2 \ell_3^2+5)-3 k^2 \ell_3^6+\ell_3^4 )   \right) \cos  (\sqrt{3}-3) ka  \\ \nonumber
    &    +k^2  \left(k^2 \ell_4^4 (k^2 \ell_3^2-1)^3-2 \ell_4^2 (k^2 \ell_3^2-1)^3+\ell_3^2 (k^4 \ell_3^4-3 k^2 \ell_3^2+3)\right)  \cos  2 (\sqrt{3}-3) ka  \Bigg]. \nonumber &&
\end{flalign}

\begin{flalign}	\label{h,minusR}
	h(k)	&=    -16 k^4 \ell_3^4 \left(- (k^2 \ell_4^2-1 )  (k^2 \ell_3^2-1 ) \cos 2 ka+k^4 \ell_4^2 \ell_3^2+k^2  (\ell_4^2-4 \ell_4\ell_3 +\ell_3^2 )+1\right)^2 . &&
\end{flalign}


\begin{flalign}	\label{w,minusR}
	w(k)	&=  4 k^2 \ell_3^2\times\Bigg[  (k^2 \ell_3^2-1 )^2 \left(2 k^2  (\ell_4^2  (k^2  (\ell_4^2  (k^2 \ell_3^2+3 )-10 \ell_3^2 )+2 )+\ell_3^2 )+6\right)  \\ \nonumber
	&  + (k^2 \ell_3^2-1 )^2 \left( (k^2  (\ell_4^2  (k^2 \ell_4^2-2 )  (k^2 \ell_3^2-1 )+\ell_3^2 )-1 )\cos  2  (\sqrt{3}-3 ) ka    -32k^2   \ell_4\ell_3   (k^2 \ell_4^2-1 )   \cos   (\sqrt{3}-3 ) ka  \right)\\ \nonumber
	&   +  (k^2 \ell_3^2-1 ) \left((k^2 \ell_4^2-1 )^2  (3 k^2 \ell_3^2+1 )^2 \cos  2  (\sqrt{3}+1 ) ka-4  (k^4 \ell_4^4-1 )   (3 k^4 \ell_3^4-2 k^2 \ell_3^2-1 ) \cos 2 \sqrt{3} ka  \right)   \\ \nonumber
	&  + 4(k^2 \ell_3^2-1 )^2 \left( ( k^2 (\ell_3^2  (k^4 \ell_4^4-1 )-k^2 \ell_4^4)+1 )\cos  2  (\sqrt{3}-2 ) ka  -6k^2  \ell_3^2  (k^2 \ell_4^2-1 )^2 \cos   (\sqrt{3}-1 ) ka \right)    \\ \nonumber
	&-4 k^2 \ell_3  (3 k^2 \ell_3^2+1 )  \left((k^2 \ell_3^3+3 \ell_3)  (k^4 \ell_4^4+6 k^2 \ell_4^2+1 )-4  (k^2 \ell_4^3+\ell_4 )(1+3 k^2  \ell_3^2 )\right)\cos  (\sqrt{3}+3 ) ka \\ \nonumber
	&  +32 k^2 \ell_3  (k^2 \ell_4^2-1 )  (k^2 \ell_3^2-1 ) \left(-m  (3 k^2 \ell_3^2+1 )+k^2 \ell_3^3+\ell_4^2  (k^4 \ell_3^3+k^2 \ell_3 )+\ell_3\right)   \cos (\sqrt{3}+1 ) ka   \\ \nonumber
	&    + (k^2 \ell_3^2-1 ) \left(8  (k^4 \ell_4^4-1 )  (k^4 \ell_3^4-1 ) \cos 2 ka-2 (k^2 \ell_4^2-1)^2 (5 k^4 \ell_3^4+2 k^2 \ell_3^2+1) \cos 4 ka\right)  \\ \nonumber
	&   -2  (k^2 \ell_3^2-1 )^2 \left(3+ k^2 (\ell_4^2  (k^2  (\ell_3^2  (k^2 \ell_4^2-10 )+3 \ell_4^2 )+2 )+\ell_3^2 )\right) \cos  2  (\sqrt{3}-1 ) ka  \\ \nonumber
	&  +4 k^2 \ell_3  (k^2 \ell_3^2-1 )^2  \left(\ell_4  (k^2 \ell_4 (\ell_4  (k^2 \ell_4\ell_3 -4 )+6 \ell_3 )-4 )+s\right) \cos  (\sqrt{3}-5) ka   \Bigg]. \nonumber &&
\end{flalign}

\subsection*{Data availability statement}

Data are available in the article.

\subsection*{Conflict of interest}

The authors have no conflict of interest.


\subsection*{Acknowledgments}
The work of M.B. was supported by the Czech Science Foundation within the project 22-18739S. P.E. was supported by the Czech Science Foundation within the project 21-07129S and by the EU project CZ.02.1.01/0.0/0.0/16\textunderscore 019/0000778. M.B. would like to thank Ji\v{r}\'{\i} Lipovsk\'{y} for helpful discussions.


\end{document}